\documentclass{article}
\usepackage{multicol,graphicx,color}
\usepackage{pslatex}
\usepackage{authblk}
\usepackage{amsthm}
\usepackage{amsmath}
\usepackage{amssymb}
\usepackage{latexsym}
\usepackage{lscape}
\usepackage{epsfig}
\usepackage{amsfonts}
\usepackage{enumerate}
\usepackage{mathrsfs}
\usepackage[
  hmarginratio={1:1},
  vmarginratio={1:1},
  textwidth=15cm,
  textheight=21cm,
  heightrounded,
]{geometry}

\usepackage[colorlinks]{hyperref}
\usepackage[numbers,sort&compress]{natbib}
\hypersetup{linkcolor=blue,citecolor=blue,filecolor=black,urlcolor=blue}

\theoremstyle{plain}
\newtheorem{theorem}{Theorem}

\newtheorem{lemma}[theorem]{Lemma}
\newtheorem{proposition}[theorem]{Proposition}
\theoremstyle{definition}

\newtheorem{remark}[theorem]{Remark}

\numberwithin{equation}{section}
\numberwithin{theorem}{section}

\newcommand{\dist}{\operatorname{dist}}
\newcommand{\supp}{\operatorname{supp}}

\newcommand{\cE}{\mathcal{E}}
\newcommand{\cX}{\mathcal{X}}
\newcommand{\cQ}{\mathcal{Q}}
\newcommand{\cZ}{\mathcal{Z}}
\author{Qun Wang\footnote{{wangq199@nenu.edu.cn}}}
\affil{School of Mathematics and Statistics, Northeast Normal University, Changchun 130024, Jilin, PR China}
\author{Xiaojun Chang\footnote{changxj100@nenu.edu.cn}}
\affil{School of Mathematics and Statistics \& Center for Mathematics and Interdisciplinary Sciences, Northeast Normal University, Changchun 130024, Jilin, PR China}
\title{Multiplicity and least energy sign-changing normalized solutions for a Kirchhoff equation on bounded domains}
\date{}

\begin{document}
\maketitle

\begin{abstract}
In this paper, we study sign-changing normalized solutions of the Kirchhoff problem
\[
\begin{cases}
-\left(a+b\displaystyle\int_\Omega|\nabla u|^2\,dx\right)\Delta u
=|u|^{p-2}u+\mu u
&\text{in }\Omega,\\
u=0
&\text{on }\partial\Omega,
\end{cases}
\]
subject to the constraint
\[
\int_\Omega |u|^2\,dx=c,
\]
where $a,b>0$, $\Omega\subset\mathbb R^N$ is a smooth bounded
connected domain, $1\le N\le3$, $c>0$ is prescribed, and
$\mu\in\mathbb R$ is a Lagrange multiplier. We assume $2<p<\infty$ for $N=1,2$ and $2<p\le6$
for $N=3$.
For every integer $k\ge1$, we obtain at least $k$ distinct pairs of sign-changing normalized solutions
for sufficiently small $c>0$. 
We determine the small mass asymptotics of their energies,
Lagrange multipliers, and rescaled profiles in terms
of the distinct Dirichlet eigenvalues and their eigenspaces.
When $2<p<2+\frac{8}{N}$, we obtain infinitely many such pairs
for every prescribed mass $c>0$, with energies tending
to $+\infty$.
At $p=2+\frac{8}{N}$, this conclusion holds for every sufficiently
small fixed mass.

We also establish the existence of a least energy
sign-changing normalized solution on general bounded domains for every prescribed mass
when $2<p<2+\frac{8}{N}$.
For sufficiently small mass, such a solution exists
on general bounded domains when $p=2+\frac{8}{N}$
and on star-shaped domains when $p>2+\frac{8}{N}$.
For any family of such least energy solutions,
the rescaled profiles approach the $L^2$-unit sphere
of the second Dirichlet eigenspace in $H_0^1(\Omega)$
as $c\to0^+$.
The finite multiplicity and least energy results include
the Sobolev critical case $N=3$, $p=6$.
The proofs combine a localized symmetric minimax
construction, energy comparisons with one-signed critical
points, and compactness arguments that preserve
the sign-changing property.
\end{abstract}

\medskip
{\small \noindent \text{Keywords:} Kirchhoff equation, normalized solution, sign-changing solution, bounded domains, variational methods.}
\medskip

\section{Introduction and main results}\label{intro}
Let $\Omega$ be a bounded connected domain in $\mathbb R^N$, $1\le N\le3$, with smooth boundary. For given $a,b>0$ and prescribed mass $c>0$, we study normalized solutions of the nonlinear Kirchhoff problem
\begin{equation}\label{K}
\begin{cases}
-\left(a+b\displaystyle\int_\Omega|\nabla u|^2\,dx\right)\Delta u=|u|^{p-2}u+\mu u
&\text{in }\Omega,\\
u=0&\text{on }\partial\Omega,\\
\displaystyle\int_\Omega|u|^2\,dx=c,
\end{cases}
\end{equation}
where $\mu\in\mathbb R$ is unknown and arises as a Lagrange multiplier. We assume that $2<p<2^*$ when $N=1,2$ and $2<p\le2^*$ when $N=3$, where
\[
2^*:=
\begin{cases}
+\infty,&N=1,2,\\
6,&N=3.
\end{cases}
\]

Normalized solutions of \eqref{K} are precisely the critical points of the energy functional
\[
E(u):=\frac{a}{2}\int_\Omega|\nabla u|^2\,dx
+\frac{b}{4}\left(\int_\Omega|\nabla u|^2\,dx\right)^2-\frac{1}{p}\int_\Omega|u|^p\,dx
\]
on the mass constraint
\[
S_c:=\left\{u\in H_0^1(\Omega):\int_\Omega|u|^2\,dx=c\right\}.
\]
The mass-critical exponent for the Kirchhoff energy
is $2+\frac{8}{N}$.
In the mass-subcritical case, $E$ is bounded from below
and coercive on $S_c$ for every $c>0$. The same conclusion holds
for sufficiently small mass in the mass-critical case.
In the mass-supercritical case, $E$ is unbounded
from below on $S_c$ for every $c>0$.

For Kirchhoff equations in the whole space $\mathbb R^N$,
normalized solutions have attracted considerable attention
\cite{Ye-2015,Li-2019,Hu-2021,Qi-2022,Cui-2022,He-2023,Feng-2023,Chen-2025}.
On bounded domains, Zhu, Li and Liang \cite{Zhu-2021}
studied a normalized Kirchhoff transmission problem
in $\mathbb R^3$.
Their results yield normalized local minimizers
for a scalar Dirichlet special case when
$\frac{14}{3}\le p\le6$ and the prescribed mass
is sufficiently small.
Wang and Chang \cite{Wang-2024} constructed positive
normalized solutions of mountain-pass type
for sufficiently small mass when
$2+\frac{8}{N}<p<2^*$.
In dimension three, Zhang \cite{Zhang-2025} obtained,
under suitable mass restrictions,
positive normalized local minimizers for
$\frac{14}{3}<p\le6$, including the Sobolev critical case,
and a second normalized solution on star-shaped domains
for $\frac{14}{3}<p<6$.

We focus on sign-changing normalized solutions of \eqref{K},
that is, solutions with $u^+\not\equiv0$ and $u^-\not\equiv0$,
where $u^+:=\max\{u,0\}$ and $u^-:=\max\{-u,0\}$.
Such a solution is called a
\emph{least energy sign-changing normalized solution}
if its energy is minimal among all sign-changing solutions
of \eqref{K} with the same prescribed mass.

For the local Schr\"odinger problem obtained by setting
$a=1$ and $b=0$ in \eqref{K}, several approaches have been
developed to obtain sign-changing normalized solutions.
Song and Zou~\cite{Song-Zou-2025,Song-Zou-2025-Correction} employed a constrained linking construction based on invariant sets of descending flow to obtain any prescribed finite number of such solutions to the Br\'ezis--Nirenberg problem in the sufficiently small mass regime. For the Sobolev subcritical pure-power problem, Song and Zou~\cite{Song-Zou-2026} combined genus theory with the method of invariant sets of descending flow. They obtained infinitely many
sign-changing normalized solutions for every mass when
$2<p<2+\frac{4}{N}$, and for sufficiently small mass when
$p=2+\frac{4}{N}$.
When $2+\frac{4}{N}<p<2^*$, they obtained any prescribed
finite number of such solutions for sufficiently small mass.
They also studied the limiting behavior of the energies
and Lagrange multipliers as the mass tends to zero.
Related linking and constrained flow techniques were
developed by Jeanjean and Song \cite{Jeanjean-Song-2025}
for $L^2$-supercritical NLS on compact metric graphs.
For further results on sign-changing normalized solutions,
we refer to \cite{Jeanjean-Lu-2019,Jeanjean-Lu-2020,Zhan-2025}
for NLS in $\mathbb R^N$,
to \cite{Xia-2023} for a Choquard equation in $\mathbb R^N$,
and to \cite{Feng-2026} for competition-diffusion
$p$-Laplacian systems on bounded domains.

An action approach was developed by De Coster, Dovetta,
Galant and Serra \cite{Coster-2024} to study least energy
normalized sign-changing solutions in the Sobolev subcritical range.
Their existence results hold for every mass when
$2<p<2+\frac{4}{N}$, for sufficiently small mass when
$p=2+\frac{4}{N}$, and for sufficiently small mass on smooth
star-shaped domains when $2+\frac{4}{N}<p<2^*$.
Galant and Weth \cite{Galant-Weth-2026} developed an
abstract framework based on ground states of
Nehari--Pankov type. Its applications include multiple
normalized solutions of small mass for equations on
bounded domains and, in the second order case,
sign-changing solutions.

For the Kirchhoff problem \eqref{K}, the positive and
negative parts are coupled through the coefficient
$a+b\int_\Omega|\nabla u|^2\,dx$.
Indeed, setting
$A_\pm:=\int_\Omega|\nabla u^\pm|^2\,dx$, we have
\[
E(u)=E(u^+)+E(u^-)+\frac{b}{2} A_+A_-.
\]
Thus, varying one nodal part also changes the diffusion
coefficient acting on the other, and the energy does not
split under separate scalings of the two parts.
 In addition, at the Sobolev critical exponent,
bounded constrained Palais--Smale sequences need not be
strongly compact. For the least energy problem, a further issue
is to prevent a minimizing sequence of sign-changing
solutions from converging to a one-signed solution.

Our first result establishes the existence and small mass
asymptotics of multiple sign-changing normalized solutions. For $2<p<2+\frac{8}{N}$, it also yields infinitely many
such solutions for every prescribed mass.
At $p=2+\frac{8}{N}$, the same conclusion holds for every
sufficiently small mass.
To state this result, let
$0<\lambda_1(\Omega)<\lambda_2(\Omega)<\cdots$
denote the increasing sequence of distinct Dirichlet
eigenvalues of $-\Delta$ in $\Omega$, and let $\cE_j$
be the eigenspace associated with $\lambda_j(\Omega)$.
For a closed subspace $X\subset H_0^1(\Omega)$, we write
\[
\mathbb S(X):=
\left\{v\in X:\int_\Omega|v|^2\,dx=1\right\}.
\]
\begin{theorem}\label{thm:multi}
Let $1\le N\le3$, $a,b>0$, and let
$\Omega\subset\mathbb R^N$ be a bounded connected domain
with a smooth boundary. Assume that one of the following
conditions holds:
\begin{enumerate}[(i)]
  \item $1\le N\le2$ and $2<p<2^*$;
  \item $N=3$ and $2<p\le6$.
\end{enumerate}
For every integer $k\ge1$, there exists $c_k^*>0$ such
that, for every $0<c<c_k^*$, problem~\eqref{K} admits
$k$ mutually distinct pairs of sign-changing
normalized solutions
$\{\pm u_{2,c}\},\ldots,\{\pm u_{k+1,c}\}$.
Moreover, for every $j\in\{2,\ldots,k+1\}$, as $c\to0^+$,
\begin{equation}\label{asym-multi}
\frac{E(u_{j,c})}{c}\to\frac{a\lambda_j(\Omega)}{2},
\qquad
\frac{\|\nabla u_{j,c}\|_2^2}{c}\to\lambda_j(\Omega),
\qquad
\mu_{j,c}\to a\lambda_j(\Omega),
\end{equation}
and
\begin{equation}\label{profile-multi}
\dist_{H_0^1}
\left(\frac{u_{j,c}}{\sqrt c},\mathbb S(\cE_j)\right)
\to0,
\end{equation}
where $\mu_{j,c}$ is the Lagrange multiplier associated
with $u_{j,c}$. 

If $2<p<2+\frac{8}{N}$, then for every fixed $c>0$,
problem~\eqref{K} admits a sequence of mutually distinct pairs of sign-changing normalized solutions
whose energies tend to $+\infty$.
If $p=2+\frac{8}{N}$, there exists $c_\infty>0$ such that
the same conclusion holds for every $0<c<c_\infty$.
\end{theorem}

Theorem~\ref{thm:multi} does not require $\Omega$ to be star-shaped, nor does it require the higher eigenvalues to be simple.  To prove the small mass results, we localize the
variational construction on subsets of $S_c$ where
$\|\nabla u\|_2^2$ is bounded by a fixed multiple of $c$.
For $c>0$ and $\rho>\lambda_1(\Omega) c$, set
\[
B_c(\rho):=
\left\{u\in S_c:\|\nabla u\|_2^2<\rho\right\},
\qquad
\overline{B_c(\rho)}:=
\left\{u\in S_c:\|\nabla u\|_2^2\le\rho\right\}.
\]
The restriction $\rho>\lambda_1(\Omega)c$ ensures that
$B_c(\rho)$ is nonempty and that the displayed closure
formula holds.
For each fixed $L>\lambda_1(\Omega)$, the Gagliardo--Nirenberg inequality
gives
\[
E(u)=\frac{a}{2}\|\nabla u\|_2^2
     +O\!\left(c^2+c^{\frac{p}{2}}\right)
\qquad\text{as }c\to0^+,
\]
uniformly for $u\in\overline{B_c(Lc)}$.
Since $p>2$, the error is of smaller order than $c$,
 which enables us to distinguish the energy levels associated with different Dirichlet eigenvalues.

We establish a symmetric minimax principle adapted to
these regions by constructing a descending deformation
that preserves the mass constraint. Using the Borsuk--Ulam theorem and spectral gap estimates,
we obtain the required intersection property and construct
energy barriers that keep the deformation within
the localized region.
The compactness needed for this construction follows from
Lemma~\ref{lem:localized-ps}.
In the Sobolev critical case, this lemma provides a
threshold $\rho^*>0$ independent of $c$.
For each fixed $L>\lambda_1(\Omega)$, the choice
$0<c<\frac{\rho^*}{L}$ gives
$\lambda_1(\Omega)c<Lc<\rho^*$,
ensuring both nonemptiness and the required compactness
on $\overline{B_c(Lc)}$.
By testing the equation with a positive first eigenfunction,
we obtain a uniform upper bound for the energies of one-signed critical points
in this region.
Since the constructed critical levels lie above this bound,
the corresponding solutions must change sign.


In the coercive regimes described above, we apply
the minimax principle on the full mass sphere to obtain
critical values tending to $+\infty$.
By testing the equation with a positive first eigenfunction
and using the Gagliardo--Nirenberg inequality, we also obtain a uniform upper bound for the energies of all one-signed constrained critical points
with the same mass.
Consequently, the critical points at sufficiently large
levels are sign-changing.

We next consider the least energy among all sign-changing
normalized solutions of a fixed mass.
By Theorem~\ref{thm:multi}, the set
\[
\mathcal K_{c}^{\mathrm{nod}}:=
\left\{
u\in S_c:
(E|_{S_c})'(u)=0,\quad
u^+\not\equiv0,\quad u^-\not\equiv0
\right\}
\]
is nonempty
for every $c>0$ when $2<p<2+\frac{8}{N}$, and for every
sufficiently small $c>0$ in the remaining admissible cases. Whenever $\mathcal K_{c}^{\mathrm{nod}}\ne\varnothing$, define
\[
d_c:=\inf_{u\in\mathcal K_{c}^{\mathrm{nod}}}E(u).
\]
Our second result establishes attainment of $d_c$ and
identifies the small mass profile of every least energy
sign-changing normalized solution.

\begin{theorem}\label{glne:main}
Let $1\le N\le3$, $a,b>0$, and let
$\Omega\subset\mathbb R^N$ be a bounded connected domain
with a smooth boundary. Assume that one of the following
conditions holds:
\begin{enumerate}[(i)]
  \item $1\le N\le3$ and $2<p\le2+\frac{8}{N}$;
  \item $1\le N\le3$, $2+\frac{8}{N}<p<2^*$, and
        $\Omega$ is star-shaped;
  \item $N=3$, $p=6$, and $\Omega$ is star-shaped.
\end{enumerate}
If $2<p<2+\frac{8}{N}$, then, for every $c>0$,
the value $d_c$ is attained by some
$w_c\in\mathcal K_{c}^{\mathrm{nod}}$.
If $p=2+\frac{8}{N}$, the same conclusion holds for every
$0<c<c_\infty$, where $c_\infty$ is the constant in
Theorem~\ref{thm:multi}.
In cases~{\rm(ii)} and~{\rm(iii)}, there exists $c^*>0$
such that the same conclusion holds for every
$0<c<c^*$.
Moreover, for any family of minimizers
$w_c\in\mathcal K_{c}^{\mathrm{nod}}$ with associated
Lagrange multipliers $\mu_c$, as $c\to0^+$,
\begin{equation}\label{asym-least-nodal}
\frac{E(w_c)}{c}\to\frac{a\lambda_2(\Omega)}{2},
\qquad
\frac{\|\nabla w_c\|_2^2}{c}\to\lambda_2(\Omega),
\qquad
\mu_c\to a\lambda_2(\Omega),
\end{equation}
and
\begin{equation}\label{profile-least-nodal}
\dist_{H_0^1}
\left(\frac{w_c}{\sqrt c},\mathbb S(\cE_2)\right)
\to0.
\end{equation}
\end{theorem}

To prove Theorem~\ref{glne:main}, we first establish
the boundedness of minimizing sequences for $d_c$.
When $2<p<2+\frac{8}{N}$, by the Gagliardo--Nirenberg
inequality it follows that $E$ is
coercive on $S_c$ for every fixed $c>0$.
The same argument applies at $p=2+\frac{8}{N}$ when
$0<c<c_\infty$. 
For $p>2+\frac{8}{N}$, we instead use the Pohozaev
identity on star-shaped domains to obtain
\[
E(u)\ge\frac{a}{4}\|\nabla u\|_2^2
\]
for constrained critical points $u\in S_c$.
In particular, at the Sobolev critical exponent,
we combine the preceding energy lower bound with
the comparison $d_c\le E(u_{2,c})=O(c)$ from
Theorem~\ref{thm:multi} to choose a minimizing sequence
satisfying $\|\nabla u_n\|_2^2=O(c)$, uniformly in $n$.
These estimates allow us to apply the localized
Palais--Smale compactness result (Lemma \ref{lem:localized-ps}), including in the
Sobolev critical case.

A further difficulty is that strong convergence alone
does not preserve the sign-changing property.
To exclude a one-signed limit, we normalize a nodal
part that could vanish along the sequence.
Using the equations and the strong convergence,
we obtain a uniform $H_0^1(\Omega)$ bound for these
normalized parts.
Together with the strong maximum principle,
this rules out one-signed limits, see
Lemma~\ref{glne:compactness}.
The limit therefore remains in
$\mathcal K_{c}^{\mathrm{nod}}$, and the continuity of $E$ yields attainment of $d_c$.

Finally, we use the rescaled equation for least energy
sign-changing normalized solutions $w_c$ of mass $c$
to establish compactness of $\frac{w_c}{\sqrt c}$
in $H_0^1(\Omega)$ as $c\to0^+$
and identify its limits as Dirichlet eigenfunctions.
Since $E(w_c)\le E(u_{2,c})$, the energy asymptotics
for $u_{2,c}$ imply that every limit of $\frac{w_c}{\sqrt c}$
has eigenvalue $\lambda_1(\Omega)$ or $\lambda_2(\Omega)$.
The same normalization argument shows that no limit belongs to the first eigenspace,
so every limit belongs to $\mathbb S(\cE_2)$.
This gives the asymptotic conclusions in
Theorem~\ref{glne:main}.

\begin{remark}\label{rem:nls}
The multiplicity and least energy arguments also apply
to the local Schr\"odinger problem obtained by setting
$a=1$ and $b=0$ in \eqref{K}, after adjusting the constants.
For the local problem, the mass-critical exponent is
$2+\frac{4}{N}$.
In the mass-supercritical regime $p>2+\frac{4}{N}$,
the least energy result holds for sufficiently small mass
provided that $\Omega$ is star-shaped. In particular, this result covers the Sobolev critical case.
The corresponding Sobolev subcritical multiplicity
and least energy existence results are contained in
\cite{Song-Zou-2026} and
\cite[Theorems~1.4 and~1.7]{Coster-2024}, respectively.
\end{remark}

\begin{remark}
In our approach, we combine the Borsuk--Ulam intersection
property with spectral estimates for small mass to
construct minimax levels on $S_c$.
The corresponding constrained critical points have
energies strictly above those of all one-signed
constrained critical points in the regions used for
the construction, and therefore change sign.
This energy comparison avoids the need to construct
invariant neighborhoods of the positive and negative
cones, which are used in the constrained linking and
genus approaches of
\cite{Song-Zou-2025,Song-Zou-2025-Correction,
Jeanjean-Song-2025,Song-Zou-2026}.

For the least energy problem, we use the critical point
equation to exclude one-signed limits of minimizing
sequences, see Lemma~\ref{glne:compactness}.
Together with localized compactness, this yields
attainment of the least energy among sign-changing
normalized solutions.

The mass remains fixed throughout these arguments.
This also distinguishes our method from the
Nehari--Pankov approach of \cite{Galant-Weth-2026},
which obtains prescribed masses from action ground
states as the frequency varies over spectral gaps.
\end{remark}

The remainder of the paper is organized as follows.
In Section~\ref{sec:multiple-solutions}, we establish
a localized symmetric minimax principle and use it
to prove Theorem~\ref{thm:multi}.
In Section~\ref{glne:section}, we prove the existence
and small mass asymptotics of least energy sign-changing
normalized solutions stated in Theorem~\ref{glne:main}.

\section{Multiple sign-changing normalized solutions}
\label{sec:multiple-solutions}

This section is devoted to the proof of
Theorem~\ref{thm:multi}. We first establish compactness
and spectral estimates for $E|_{S_c}$, then prove a
localized symmetric minimax principle, and finally
apply it to construct solutions associated with the
distinct Dirichlet eigenvalues.

\subsection{Compactness and spectral estimates}
\label{subsec:preliminaries}

For \(1\le q<\infty\), we denote the usual \(L^q(\Omega)\) norm by
\[
\|u\|_q:=\left(\int_\Omega |u|^q\,dx\right)^{\frac{1}{q}}.
\]
We equip \(H_0^1(\Omega)\) with the norm
\[
\|u\|:=\left(\int_\Omega |\nabla u|^2\,dx\right)^{\frac{1}{2}}.
\]
By the Poincar\'e inequality, this norm is equivalent to the standard \(H^1(\Omega)\) norm. For $2<q<2^*$, and also for $q=6$ when $N=3$, let $\mathcal S_q>0$ denote the optimal constant for the Sobolev embedding, namely,
\begin{equation}\label{Sobolev}
\|u\|_q
\le
\mathcal S_q\|u\|,
\qquad
u\in H_0^1(\Omega).
\end{equation}

Set
\[
q_p:=\frac{N(p-2)}{4},\qquad r:=\frac{p-2}{2}.
\]
For $2<p<2^*$, the Gagliardo--Nirenberg inequality yields a constant $C_{N,p}>0$ such that
\begin{equation}\label{GN}
\|u\|_p^p\le C_{N,p}\|u\|_2^{p-2q_p}\|\nabla u\|_2^{2q_p},\qquad u\in H_0^1(\Omega).
\end{equation}
When $N=3$ and $p=6$, we have $q_p=3$ and set
$C_{3,6}:=\mathcal S_6^6$, so that \eqref{GN}
also holds at the Sobolev critical exponent.
Consequently, in all the cases under consideration,
for every $L>\lambda_1(\Omega)$,
\begin{equation}\label{eq:small-gradient-nonlinear}
\|u\|_p^p\le C_{N,p}L^{q_p}c^{\frac{p}{2}} \qquad\text{for every }u\in\overline{B_c(Lc)}.
\end{equation}

For $u\in S_c$, we write
\[
T_uS_c:=\left\{v\in H_0^1(\Omega):\int_\Omega uv\,dx=0\right\},
\]
and denote the constrained derivative of $E$ on $S_c$ by $(E|_{S_c})'(u)$. Its norm is given by
\[
\|(E|_{S_c})'(u)\|_{T_u^*S_c}:=\sup_{\substack{v\in T_uS_c\\ \|v\|=1}}|(E|_{S_c})'(u)[v]|.
\]

Every constrained critical point $u\in S_c$ satisfies, for some
$\mu\in\mathbb R$,
\[
-M_u\Delta u=|u|^{p-2}u+\mu u
\quad\text{in }\Omega,
\qquad
u=0\quad\text{on }\partial\Omega,
\]
in the weak sense, where
$M_u:=a+b\|\nabla u\|_2^2>0$. By the Br\'ezis--Kato estimate and elliptic regularity arguments (see \cite{BK, Trudinger}), we get 
$u\in C^{2,\alpha}(\overline\Omega)$ for some $\alpha\in (0,1)$.
In particular, using the strong maximum principle, every nonnegative constrained critical point
is positive in $\Omega$.
This regularity also justifies the Pohozaev identity used below.

We first establish a Palais--Smale compactness result.
In the Sobolev critical case $N=3$, $p=6$, set
$\rho^*:=\frac{\sqrt a}{\mathcal S_6^3}$. 

\begin{lemma}\label{lem:localized-ps}
Let $c>0$, $\rho>\lambda_1(\Omega)c$, and let $\{u_n\}\subset\overline{B_c(\rho)}$ satisfy
\[
E(u_n)\to d,\qquad\|(E|_{S_c})'(u_n)\|_{T_{u_n}^*S_c}\to0
\]
for some $d\in\mathbb R$. Assume either
\begin{enumerate}[(i)]
\item $1\leq N\leq3$ and $2<p<2^*$, or
\item $N=3$, $p=6$, and $\rho\le\rho^*$.
\end{enumerate}
Then, $\{u_n\}$ has a strongly convergent subsequence in $H_0^1(\Omega)$. In particular, $E|_{S_c}$ satisfies the $(PS)_d$ condition in $\overline{B_c(\rho)}$ under the corresponding assumption.
\end{lemma}
\begin{proof}
For $2<p<2^*$, the conclusion follows from the
argument in the proof of \cite[Lemma~4]{Wang-2024},
which uses only boundedness in $H_0^1(\Omega)$
and the compact embeddings into $L^2(\Omega)$
and $L^p(\Omega)$.

Assume now that $N=3$ and $p=6$. Set $e_n:=\|\nabla u_n\|_2^2$. Since $\{u_n\}\subset\overline{B_c(\rho)}$ is bounded in $H_0^1(\Omega)$, we may assume, up to a subsequence, that $e_n\to A$ for some $A\in[0,\rho]$ and that there exists $u\in S_c$ such that
\begin{align*}
&u_n\rightharpoonup u \ \text{ in }H_0^1(\Omega),\\
&u_n\to u \ \text{ in }L^q(\Omega),\qquad 2\le q<6,\\
&u_n(x)\to u(x) \ \text{ for a.e. }x\in\Omega.
\end{align*}
Define the linear functional $\Psi(u):H_0^1(\Omega)\to\mathbb R$ by
\[
\Psi(u)v
:=
a\int_\Omega\nabla u\cdot\nabla v\,dx
+b\|\nabla u\|_2^2
 \int_\Omega\nabla u\cdot\nabla v\,dx
-\int_\Omega |u|^{p-2}uv\,dx .
\]
Since $\{u_n\}$ is a Palais--Smale sequence for $E|_{S_c}$, there exists a sequence $\{\mu_n\}\subset\mathbb R$ such that
\begin{equation}\label{eq:ps-multiplier}
\Psi(u_n)v-\mu_n\int_\Omega u_nv\,dx=o_n(1)\|v\|,
\qquad v\in H_0^1(\Omega),
\end{equation}
where $o_n(1)\to0$ uniformly for $\|v\|\le1$. Taking $v=u_n$ in \eqref{eq:ps-multiplier}, we obtain
\[
\mu_n c
=
(a+b e_n)e_n-\|u_n\|_p^p+o_n(1).
\]
Hence $\{\mu_n\}$ is bounded, and, up to a subsequence, $\mu_n\to\mu$ for some $\mu\in\mathbb R$.

We now apply a splitting argument based on the Br\'ezis--Lieb lemma \cite{Brezis-Lieb-1983}. Set $v_n:=u_n-u$. Then $v_n\rightharpoonup0$ in $H_0^1(\Omega)$, $v_n\to0$ in $L^2(\Omega)$ and a.e. in $\Omega$. Moreover,
\[
e_n
=
\|\nabla u\|_2^2+\|\nabla v_n\|_2^2+o(1).
\]
By the Br\'ezis--Lieb lemma,
\[
\|u_n\|_{6}^{6}=\|u\|_{6}^{6}+\|v_n\|_{6}^{6}+o(1).
\]
Since $\{u_n\}$ is bounded in $L^6(\Omega)$ and $u_n\to u$ almost everywhere, we also have $|u_n|^4u_n\rightharpoonup |u|^4u$ in $L^{\frac{6}{5}}(\Omega)$. Hence,
\[
\int_\Omega |u_n|^{4}u_n v_n\,dx=\|v_n\|_{6}^{6}+o(1).
\]
Taking $v=v_n$ in \eqref{eq:ps-multiplier}, and using $e_n\to A$ and $v_n\to0$ in $L^2(\Omega)$, we obtain
\begin{equation}\label{eq:critical-remainder}
(a+bA)\|\nabla v_n\|_2^2=\|v_n\|_{6}^{6}+o(1).
\end{equation}

We distinguish two cases. Either $\|v_n\|_{6}^{6}\to0$, or, up to a subsequence, $\|v_n\|_{6}^{6}\to D>0$. In the first case, using \eqref{eq:critical-remainder} and $a+bA\ge a>0$, we have $\|\nabla v_n\|_2\to0$, and hence $u_n\to u$ in $H_0^1(\Omega)$. Assume now that the second case occurs. By \eqref{eq:critical-remainder}, up to a subsequence,
\[
\|\nabla v_n\|_2^2\to\ell>0,
\qquad
D=(a+bA)\ell.
\]
From the Sobolev inequality \eqref{Sobolev}, it follows that
\[
(a+bA)\ell=D\le\mathcal S_6^6\ell^3,
\]
which implies
\[
\ell\ge\frac{\sqrt{a+bA}}{\mathcal S_6^3}\ge\rho^*.
\]
On the other hand, since $u\in S_c$, the
Poincar\'e inequality yields
\[
\|\nabla u\|_2^2\ge\lambda_1(\Omega)c>0.
\]
Therefore,
\[
\ell
=A-\|\nabla u\|_2^2
\le\rho-\lambda_1(\Omega)c
\le\rho^*-\lambda_1(\Omega)c
<\rho^*,
\]
contradicting $\ell\ge\rho^*$.
Thus, the second case cannot occur, and
$u_n\to u$ in $H_0^1(\Omega)$.
\end{proof}

We next derive an energy estimate for one-signed constrained critical points. This estimate separates one-signed critical points from the higher minimax levels. For $L>\lambda_1(\Omega)$, set
\[
c_1(L):=\min\left\{1,
\left(\frac{a\bigl(\lambda_2(\Omega)-\lambda_1(\Omega)\bigr)}{4C_{N,p}L^{q_p}}\right)^{\frac{1}{r}},
\frac{a\bigl(\lambda_2(\Omega)-\lambda_1(\Omega)\bigr)}{2bL^2}\right\}.
\]

\begin{lemma}\label{lem:one-signed-energy}
Let $L>\lambda_1(\Omega)$ and $0<c<c_1(L)$. Then, for every one-signed constrained critical point $u\in\overline{B_c(Lc)}$, one has
\[
E(u)<\frac{a}{4}\bigl(\lambda_1(\Omega)+\lambda_2(\Omega)\bigr)c.
\]
\end{lemma}

\begin{proof}
Let $u\in\overline{B_c(Lc)}$ be a one-signed constrained critical point. By the Lagrange multiplier rule, there exists $\mu\in\mathbb R$ such that
\begin{equation}\label{eq:solution}
-\left(a+b\|\nabla u\|_2^2\right)\Delta u=|u|^{p-2}u+\mu u\qquad\text{in }\Omega.
\end{equation}
After replacing $u$ by $-u$ if necessary, we may assume that $u\ge0$. From elliptic regularity and the strong maximum principle, we have $u>0$ in $\Omega$. Testing \eqref{eq:solution} with the positive $L^2$-normalized first eigenfunction $\phi_1$, it follows that
\[
(a+b\|\nabla u\|_2^2)\int_\Omega \nabla u\cdot\nabla\phi_1\,dx=\int_\Omega u^{p-1}\phi_1\,dx+\mu\int_\Omega u\phi_1\,dx.
\]
Since
\[
\int_\Omega \nabla u\cdot\nabla\phi_1\,dx=\lambda_1(\Omega)\int_\Omega u\phi_1\,dx,
\]
we obtain
\[
\mu=(a+b\|\nabla u\|_2^2)\lambda_1(\Omega)
-\frac{\int_\Omega u^{p-1}\phi_1\,dx}{\int_\Omega u\phi_1\,dx}<(a+b\|\nabla u\|_2^2)\lambda_1(\Omega).
\]
Testing~\eqref{eq:solution} with $u$ and using \eqref{eq:small-gradient-nonlinear}, we have
\[
\|\nabla u\|_2^2<\lambda_1(\Omega)c+\frac{1}{a}\|u\|_p^p\le
\lambda_1(\Omega)c+\frac{C_{N,p}}{a}L^{q_p}c^{\frac{p}{2}}.
\]
Consequently, for $0<c<c_1(L)$,
\[
E(u)<\frac{a\lambda_1(\Omega)}{2}c+\frac{C_{N,p}}{2}L^{q_p}c^{\frac{p}{2}}
+\frac{b}{4}L^2c^2<\frac{a}{4}(\lambda_2(\Omega)+\lambda_1(\Omega))c.
\]
\end{proof}

The next lemma characterizes the asymptotic behavior of constrained critical points as $c\to0^+$ in terms of the Dirichlet spectrum.

\begin{lemma}\label{lem:spectral-identification}
Let $j\ge1$, $L>\lambda_j(\Omega)$, and $\bar c>0$. For every $c\in(0,\bar c)$, let $u_c\in\overline{B_c(Lc)}$ be a constrained critical point of \(E|_{S_c}\), with corresponding Lagrange multiplier \(\mu_c\), namely,
\begin{equation}\label{eq:spectral-family-equation}
-\left(a+b\|\nabla u_c\|_2^2\right)\Delta u_c=|u_c|^{p-2}u_c+\mu_cu_c\qquad\text{in }\Omega.
\end{equation}
Assume that
\begin{equation}\label{eq:spectral-family-energy}
\frac{E(u_c)}{c}\to\frac{a\lambda_j(\Omega)}{2}\qquad\text{as }c\to0^+.
\end{equation}
Then, as $c\to0^+$,
\begin{equation*}
\frac{\|\nabla u_c\|_2^2}{c}\to\lambda_j(\Omega),
\qquad\mu_c\to a\lambda_j(\Omega),\qquad\operatorname{dist}_{H_0^1}
\left(\frac{u_c}{\sqrt c},\mathbb S(\mathcal E_j)\right)\to0.
\end{equation*}
\end{lemma}
\begin{proof}
Set $v_c:=\frac{u_c}{\sqrt c}$. Then $\|v_c\|_2=1$, and from $u_c\in\overline{B_c(Lc)}$ we obtain $\|\nabla v_c\|_2^2=\frac{\|\nabla u_c\|_2^2}{c}\le L$. Hence $\{v_c\}_{0<c<\bar c}$ is uniformly bounded in $H_0^1(\Omega)$ and therefore also in $L^p(\Omega)$. Moreover,
\[
\frac{E(u_c)}{c}=\frac{a}{2}\|\nabla v_c\|_2^2+\frac{b}{4}c\|\nabla v_c\|_2^4
-\frac{1}{p}c^{\frac{p-2}{2}}\|v_c\|_p^p.
\]
Since $p>2$, the last two terms tend to zero as $c\to0^+$. Together with \eqref{eq:spectral-family-energy}, we have $\|\nabla v_c\|_2^2\to\lambda_j(\Omega)$, so $\frac{\|\nabla u_c\|_2^2}{c}\to\lambda_j(\Omega)$.

Dividing \eqref{eq:spectral-family-equation} by $\sqrt c$, we can see that
\begin{equation}\label{eq:normalized-spectral-family-equation}
-\left(a+bc\|\nabla v_c\|_2^2\right)\Delta v_c=c^{\frac{p-2}{2}}|v_c|^{p-2}v_c+\mu_cv_c.
\end{equation}
Testing \eqref{eq:normalized-spectral-family-equation} with $v_c$, we get
\[
\mu_c=\left(a+bc\|\nabla v_c\|_2^2\right)\|\nabla v_c\|_2^2-c^{\frac{p-2}{2}}\|v_c\|_p^p.
\]
Hence $\mu_c\to a\lambda_j(\Omega)$.

It remains to identify the limiting profile. Let $c_n\to0^+$ be arbitrary and set $v_n:=v_{c_n}$. Since $\{v_n\}$ is bounded in $H_0^1(\Omega)$, after passing to a subsequence there exists $v\in H_0^1(\Omega)$ such that
\[
v_n\rightharpoonup v\quad\text{in }H_0^1(\Omega),
\qquad v_n\to v \quad\text{in }L^q(\Omega) \quad\text{for every }2\le q<2^*.
\]
In particular, $\|v\|_2=1$. Since $\{v_n\}$ is bounded in $L^p(\Omega)$ and $p>2$, we have $c_n^{\frac{p-2}{2}}|v_n|^{p-2}v_n\to0$ in $H^{-1}(\Omega)$. Passing to the limit in \eqref{eq:normalized-spectral-family-equation}, we obtain
\[
-\Delta v=\lambda_j(\Omega)v\qquad\text{in }\Omega.
\]
Thus $v\in\mathcal E_j$, and since $\|v\|_2=1$, we have $v\in\mathbb S(\mathcal E_j)$. Moreover, $\|\nabla v_n\|_2^2\to\lambda_j(\Omega)=\|\nabla v\|_2^2$. Together with $v_n\rightharpoonup v$ in $H_0^1(\Omega)$, this implies $v_n\to v$ in $H_0^1(\Omega)$.
Since $c_n\to0^+$ was arbitrary, it follows that
\[
\operatorname{dist}_{H_0^1}\left(\frac{u_c}{\sqrt c},\mathbb S(\mathcal E_j)\right)
\to0\qquad\text{as }c\to0^+.
\]
\end{proof}

\subsection{A localized symmetric minimax principle}
\label{subsec:abstract-minimax}

We work in the following standard constrained variational framework (see \cite{Berestycki-1983}). Let $X$ be a real Banach space and let $H$ be a real Hilbert space. We assume that
\[
X\hookrightarrow H\hookrightarrow X^*
\]
continuously, where $H$ is identified with its dual. We denote the norms of $X$ and $H$ by $\|\cdot\|_X$ and $\|\cdot\|_H$, respectively, and the inner product of $H$ by $(\cdot,\cdot)_H$. For $c>0$, we consider the constraint manifold
\[
M_c:=\{u\in X:\|u\|_H^2=c\}.
\]
We equip $M_c$ with the topology and distance induced
by $\|\cdot\|_X$.
It is a smooth submanifold of $X$ of codimension one,
with tangent space
\[
T_uM_c=\{v\in X:(u,v)_H=0\},
\qquad u\in M_c.
\]
In particular, $T_{-u}M_c=T_uM_c$.

Let $J\in C^1(X,\mathbb R)$ be even, and let $J|_{M_c}$ denote its restriction to $M_c$. Then $J|_{M_c}$ is a $C^1$ functional on $M_c$, and for any $u\in M_c$,
\[
\langle (J|_{M_c})'(u),v\rangle
=
\langle J'(u),v\rangle,
\qquad v\in T_uM_c.
\]
The norm of the constrained derivative is given by
\[
\|(J|_{M_c})'(u)\|_{T_u^*M_c}
:=
\sup_{\substack{v\in T_uM_c\\ \|v\|_X=1}}
|\langle J'(u),v\rangle|.
\]
We call $u\in M_c$ a constrained critical point of $J$ if $(J|_{M_c})'(u)=0$, and set
\[
\mathcal K_c:=\{u\in M_c:(J|_{M_c})'(u)=0\},\qquad\widetilde M_c:=M_c\setminus\mathcal K_c.
\]
For $\zeta\in\mathbb R$, we also set
\[
\mathcal K_{c,\zeta}:=\{u\in\mathcal K_c:J(u)=\zeta\},\qquad J^\zeta:=\{u\in M_c:J(u)\le\zeta\}.
\]
Let $\mathcal C$ be a closed subset of $M_c$. We say that $J|_{M_c}$ satisfies the $(PS)_\zeta$ condition in $\mathcal C$ if every sequence $\{u_n\}\subset\mathcal C$ such that
\[
J(u_n)\to\zeta,\qquad\|(J|_{M_c})'(u_n)\|_{T_{u_n}^*M_c}\to0
\]
has a subsequence converging strongly in $X$ to a point of $\mathcal C$.

We first construct an odd pseudo-gradient vector field tangent to $M_c$.
\begin{lemma}
\label{lem:abstract-odd-pseudogradient}
There exists a locally Lipschitz vector field $Y:\widetilde M_c\to X$ such that
\[
Y(u)\in T_uM_c,\qquad Y(-u)=-Y(u),
\]
and
\begin{equation}\label{eq:abstract-pseudogradient}
\|Y(u)\|_X\le2\|(J|_{M_c})'(u)\|_{T_u^*M_c},
\qquad (J|_{M_c})'(u)[Y(u)]\ge\frac{1}{2}\|(J|_{M_c})'(u)\|_{T_u^*M_c}^2
\end{equation}
for every $u\in\widetilde M_c$.
\end{lemma}

\begin{proof}
We adapt the standard pseudo-gradient construction
(see \cite{Berestycki-1983,Willem-1996}) to $J|_{M_c}$.
For $u\in M_c$, set
\[
g(u):=\|(J|_{M_c})'(u)\|_{T_u^*M_c},
\qquad
P_uv:=v-\frac{(u,v)_H}{c}\,u,
\quad v\in X.
\]
Then $P_u$ is a bounded projection onto $T_uM_c$,
and the map $u\mapsto P_u$ is locally Lipschitz
from $M_c$ into $\mathcal L(X)$.
The continuity of $J'$ and of these projections
implies that $g$ is continuous on $M_c$.
Moreover, $g>0$ on $\widetilde M_c$.

For each $w\in\widetilde M_c$, choose $e_w\in T_wM_c$
such that $\|e_w\|_X=1$ and
$J'(w)[e_w]>\frac{3}{4}g(w)$.
The field $V_w(u):=g(w)P_ue_w$ is locally Lipschitz
and tangent to $M_c$.
At $u=w$, we have $P_we_w=e_w$, so
\[
\|V_w(w)\|_X=g(w),
\qquad
J'(w)[V_w(w)]>\frac{3}{4}g(w)^2.
\]
By continuity, there is an open neighborhood $U_w$
of $w$ in $\widetilde M_c$ such that
\[
\|V_w(u)\|_X\le2g(u),
\qquad
J'(u)[V_w(u)]\ge\frac{1}{2}g(u)^2,
\qquad u\in U_w.
\]
\begin{sloppypar}
Using a locally finite, locally Lipschitz partition
of unity subordinate to $\{U_w\}_{w\in\widetilde M_c}$,
we combine these local fields into a locally Lipschitz
tangent field $Y_0:\widetilde M_c\to X$ satisfying
\end{sloppypar}
\[
\|Y_0(u)\|_X\le2g(u),
\qquad
J'(u)[Y_0(u)]\ge\frac{1}{2}g(u)^2,
\qquad u\in\widetilde M_c.
\]

The evenness of $J$ implies that $\widetilde M_c$
is symmetric, $J'(-u)=-J'(u)$, and $g(-u)=g(u)$.
We can therefore define
\[
Y(u):=\frac{1}{2}\bigl(Y_0(u)-Y_0(-u)\bigr).
\]
The field $Y$ is odd and locally Lipschitz.
It is also tangent to $M_c$, because
$T_{-u}M_c=T_uM_c$. 
The norm estimate in \eqref{eq:abstract-pseudogradient}
follows from the triangle inequality, while
\[
J'(u)[Y(u)]
=\frac{1}{2}\bigl(
J'(u)[Y_0(u)]+J'(-u)[Y_0(-u)]
\bigr)
\ge\frac{1}{2}g(u)^2.
\]
Since $Y(u)\in T_uM_c$, this proves
\eqref{eq:abstract-pseudogradient}.
\end{proof}

We now use this vector field to establish a deformation lemma on $M_c$.

\begin{lemma}
\label{lem:abstract-deformation}
Let $\mathcal C\subset\mathcal G\subset M_c$, where $\mathcal C$ is symmetric and closed in $M_c$ and $\mathcal G$ is symmetric and open in $M_c$.  Let $\zeta\in\mathbb R$. Assume that
\begin{enumerate}[(i)]
\item $J|_{M_c}$ satisfies the $(PS)_\zeta$ condition in $\mathcal C$.
\item $\mathcal K_{c,\zeta}\cap\mathcal C=\varnothing$.
\item for some $\delta>0$,
\begin{equation}\label{eq:abstract-deformation-wall}
 \zeta+4\delta<\inf_{u\in\overline{\mathcal G\setminus\mathcal C}}J(u),
\end{equation}
where the closure is taken in $M_c$.
\end{enumerate}
Then, there exist $\varepsilon\in(0,\delta)$, a time $T>0$, and a continuous flow $\eta:[0,\infty)\times M_c\to M_c$ such that
\begin{enumerate}[(i)]
\item $\eta(0,u)=u$ and $\eta(t,-u)=-\eta(t,u)$.
\item $J(\eta(t,u))\le J(u)$ for every $t\ge0$.
\item $\eta\bigl(t,\mathcal C\cap J^{\zeta+\varepsilon}\bigr)\subset\mathcal C\cap J^{\zeta+\varepsilon}$ for every $t\ge0$.
\item $\eta\bigl(T,\mathcal C\cap J^{\zeta+\varepsilon}\bigr)\subset\mathcal C\cap J^{\zeta-\varepsilon}$.
\end{enumerate}
\end{lemma}

\begin{proof}
By the $(PS)_\zeta$ condition and the assumption $\mathcal K_{c,\zeta}\cap\mathcal C=\varnothing$, there exist $\widehat\varepsilon>0$ and $\sigma_*>0$ such that $\|(J|_{M_c})'(u)\|_{T_u^*M_c}\ge\sigma_*$ whenever $u\in\mathcal C$ and $|J(u)-\zeta|\le2\widehat\varepsilon$. Choose $0<\varepsilon<\min\{\widehat\varepsilon,\delta\}$. If $\mathcal C\cap J^{-1}([\zeta-\varepsilon,\zeta+\varepsilon])=\varnothing$, then $\mathcal C\cap J^{\zeta+\varepsilon}\subset J^{\zeta-\varepsilon}$. Thus, $\eta(t,u)=u$ for all $t\ge0$, with any $T>0$,
satisfies {\rm(i)--(iv)}. We may therefore assume that $\mathcal C\cap J^{-1}([\zeta-\varepsilon,\zeta+\varepsilon])\ne\varnothing$.

Let
\[
\begin{aligned}
 U_1&:=\left\{u\in M_c:\|(J|_{M_c})'(u)\|_{T_u^*M_c}>\frac{3\sigma_*}{4}\right\},\\
 U_2&:=\left\{u\in M_c:\|(J|_{M_c})'(u)\|_{T_u^*M_c}>\frac{\sigma_*}{2}\right\}.
\end{aligned}
\]
Then $\mathcal C\cap J^{-1}([\zeta-\varepsilon,\zeta+\varepsilon])\subset U_1$ and $\overline{U_1}\subset U_2\subset\widetilde M_c$. If $U_1=M_c$, we set $\vartheta\equiv1$. Otherwise, we define
\[
 \vartheta(u):=
\frac{\dist(u,M_c\setminus U_1)}{\dist(u,\mathcal C\cap J^{-1}([\zeta-\varepsilon,\zeta+\varepsilon]))+\dist(u,M_c\setminus U_1)}.
\]
Hence, $\vartheta:M_c\to[0,1]$ is locally Lipschitz and even, $\vartheta=1$ on $\mathcal C\cap J^{-1}([\zeta-\varepsilon,\zeta+\varepsilon])$, and $\supp\vartheta\subset\overline{U_1}\subset U_2$.  Choose $\varphi\in C^1(\mathbb R,[0,1])$ such that
\[
 \varphi(s)=1\quad\text{if }|s-\zeta|\le\varepsilon,
 \qquad\varphi(s)=0\quad\text{if }|s-\zeta|\ge2\varepsilon.
\]
Let $Y$ be the odd pseudo-gradient vector field given by Lemma~\ref{lem:abstract-odd-pseudogradient}. Define
\[
\widetilde Y(u):=\vartheta(u)\varphi(J(u))\frac{Y(u)}{1+\|Y(u)\|_X},\qquad u\in U_2,
\]
and set $\widetilde Y(u):=0$ for $u\notin U_2$. Since $\supp\vartheta\subset U_2$, the vector field $\widetilde Y$ is globally bounded, locally Lipschitz, odd, and tangent to $M_c$, with $\|\widetilde Y(u)\|_X\le1$.

We consider the Cauchy problem
\begin{equation}\label{eq:abstract-flow-equation}
\frac{d}{dt}\eta(t,u)=-\widetilde Y(\eta(t,u)),\qquad\eta(0,u)=u.
\end{equation}
Since $\widetilde Y$ is locally Lipschitz, it generates a unique local flow on $M_c$ which depends continuously on the initial value. Moreover, the tangency of $\widetilde Y$ to $M_c$ implies that
\[
\frac{d}{dt}\|\eta(t,u)\|_H^2=-2\bigl(\eta(t,u),\widetilde Y(\eta(t,u))\bigr)_H=0,
\]
and hence the flow preserves $M_c$. Since $\widetilde Y$ is bounded and $M_c$ is closed in $X$, the usual continuation argument shows that the solution is defined for every $t\ge0$. By the oddness of $\widetilde Y$ and the uniqueness of solutions to \eqref{eq:abstract-flow-equation}, we have $\eta(t,-u)=-\eta(t,u)$. Moreover,
\begin{equation}\label{eq:abstract-energy-monotonicity}
\frac{d}{dt}J(\eta(t,u))=-(J|_{M_c})'(\eta(t,u))\bigl[\widetilde Y(\eta(t,u))\bigr]\le0.
\end{equation}

We next prove the invariance of $\mathcal C\cap J^{\zeta+\varepsilon}$. Let $u\in\mathcal C\cap J^{\zeta+\varepsilon}$. By \eqref{eq:abstract-energy-monotonicity}, $J(\eta(t,u))\le\zeta+\varepsilon<\zeta+4\delta$ for every $t\ge0$. Since $\mathcal C\subset\mathcal G$ and $\mathcal G$ is open in $M_c$, we have $\partial_{M_c}\mathcal C\subset\overline{\mathcal G\setminus\mathcal C}$. Hence, \eqref{eq:abstract-deformation-wall} implies that $\eta(t,u)$ remains in $\mathcal C$ for every $t\ge0$. Therefore,
\[
\eta\bigl(t,\mathcal C\cap J^{\zeta+\varepsilon}\bigr)\subset
\mathcal C\cap J^{\zeta+\varepsilon}\qquad\text{for every }t\ge0,
\]
which proves {\rm(iii)}.

To prove {\rm(iv)}, fix $u\in\mathcal C\cap J^{\zeta+\varepsilon}$. By {\rm(iii)} and \eqref{eq:abstract-energy-monotonicity}, as long as $J(\eta(t,u))>\zeta-\varepsilon$, we have
\[
\eta(t,u)\in\mathcal C\cap J^{-1}([\zeta-\varepsilon,\zeta+\varepsilon]).
\]
Consequently, $\vartheta(\eta(t,u))=\varphi(J(\eta(t,u)))=1$. Using \eqref{eq:abstract-pseudogradient} and the lower bound for $\|(J|_{M_c})'(u)\|_{T_u^*M_c}$, we obtain
\begin{equation}\label{eq:uniform-energy-descent}
\frac{d}{dt}J(\eta(t,u))\le-\frac{1}{2}\frac{\|(J|_{M_c})'(\eta(t,u))\|_{T_{\eta(t,u)}^*M_c}^2}{1+2\|(J|_{M_c})'(\eta(t,u))\|_{T_{\eta(t,u)}^*M_c}}\le-\kappa_*,
\end{equation}
where
\[
\kappa_*:=\frac{\sigma_*^2}{2(1+2\sigma_*)}>0.
\]
We set $T:=\frac{2\varepsilon}{\kappa_*}$. If $J(\eta(T,u))>\zeta-\varepsilon$, then \eqref{eq:uniform-energy-descent} holds for every $t\in[0,T]$. Therefore
\[
J(\eta(T,u))\le J(u)-\kappa_* T\le\zeta+\varepsilon-2\varepsilon=\zeta-\varepsilon,
\]
a contradiction. Thus $J(\eta(T,u))\le\zeta-\varepsilon$, and consequently
\[
\eta\bigl(T,\mathcal C\cap J^{\zeta+\varepsilon}\bigr)\subset\mathcal C\cap J^{\zeta-\varepsilon}.
\]
The proof is complete.
\end{proof}

We now derive a symmetric minimax principle from
Lemma~\ref{lem:abstract-deformation}. Let $Q\subset M_c$ be compact and symmetric, and let $\mathcal C\subset\mathcal G\subset M_c$ be symmetric sets, with $\mathcal C$ closed in $M_c$ and $\mathcal G$ open in $M_c$. Fix a subset $\mathcal Z\subset\mathcal G$ and define the class of odd continuous maps
\[
\Gamma:=\{h\in C(Q,\mathcal G):h(-q)=-h(q)\text{ for all }q\in Q\}.
\]
Assume that
\begin{enumerate}
\renewcommand{\labelenumi}{(A\arabic{enumi})}
\item $\Gamma\ne\varnothing$, and every $h\in\Gamma$ satisfies $h(Q)\cap\mathcal Z\ne\varnothing$.

\item There exist $\beta\le\alpha$ and $h_0\in\Gamma$ such that
\[
\inf\limits_{u\in\mathcal Z}J(u)\ge\beta
\quad\text{and}\quad
\max\limits_{q\in Q}J(h_0(q))\le\alpha.
\]

\item There exists $\delta>0$ such that
$\alpha+4\delta<\inf\limits_{u\in\overline{\mathcal G\setminus\mathcal C}}J(u)$.

\item For every $\zeta\in[\beta,\alpha]$, $J|_{M_c}$ satisfies the $(PS)_\zeta$ condition in
$\mathcal C$.
\end{enumerate}

We introduce the following minimax level:
\begin{equation*}
 \mathfrak m:=\inf_{h\in\Gamma}\max_{q\in Q}J(h(q)).
\end{equation*}

\begin{proposition}
\label{prop:localized-minimax}
Assume \textup{(A1)--(A4)}.
Then $\beta\le\mathfrak m\le\alpha$, and there exists
$u_*\in\mathcal C$ such that 
\[
 J(u_*)=\mathfrak m,\qquad (J|_{M_c})'(u_*)=0.
\]
\end{proposition}

\begin{proof}
By ${\rm(A1)}$ and ${\rm(A2)}$, every $h\in\Gamma$ satisfies
\[
 \max_{q\in Q}J(h(q))\ge\inf_{u\in\mathcal Z}J(u)\ge\beta.
\]
Using the test map $h_0$, we obtain $\beta\le\mathfrak m\le\alpha$.

Suppose, by contradiction, that $\mathcal K_{c,\mathfrak m}\cap\mathcal C=\varnothing$. Since $\mathfrak m\in[\beta,\alpha]$, assumptions {\rm(A3)} and {\rm(A4)}, together with $\mathcal K_{c,\mathfrak m}\cap\mathcal C=\varnothing$, allow us to apply Lemma~\ref{lem:abstract-deformation} with $\zeta=\mathfrak m$. We obtain $\varepsilon\in(0,\delta)$, $T>0$, and an odd descending flow $\eta:[0,\infty)\times M_c\to M_c$. By the definition of $\mathfrak m$, we may choose $h\in\Gamma$ such that
\[
 \max_{q\in Q}J(h(q))<\mathfrak m+\varepsilon.
\]
We claim that $h(Q)\subset\mathcal C$. Indeed, if $h(q_0)\in\mathcal G\setminus\mathcal C$ for some $q_0\in Q$, then it follows from {\rm(A3)} that
\[
 J(h(q_0))\ge\inf_{u\in\overline{\mathcal G\setminus\mathcal C}}J(u)>\alpha+4\delta\ge\mathfrak m+4\delta,
\]
which contradicts the choice of $h$. Hence $h(Q)\subset\mathcal C\cap J^{\mathfrak m+\varepsilon}$. Define $\widetilde h(q):=\eta(T,h(q))$. By Lemma~\ref{lem:abstract-deformation}{\rm(iv)}, $\widetilde h(Q)\subset\mathcal C\subset\mathcal G$. Since both $h$ and $\eta(T,\cdot)$ are odd, we have $\widetilde h\in\Gamma$. Therefore,
\[
 \max_{q\in Q}J(\widetilde h(q))\le\mathfrak m-\varepsilon,
\]
which contradicts the definition of $\mathfrak m$. Thus $\mathcal K_{c,\mathfrak m}\cap\mathcal C\ne\varnothing$, and the conclusion follows.
\end{proof}


\subsection{Spectral minimax levels}
\label{subsec:spectral-constructions}

We now apply Proposition~\ref{prop:localized-minimax} to the Kirchhoff functional and then prove that the resulting critical points are sign-changing by a separate energy comparison. Fix $j\ge2$ and choose $\lambda_j(\Omega)<\Lambda_j<\widehat{\Lambda}_j<\lambda_{j+1}(\Omega)$. Set
\[
\mathcal X_j:=\mathcal E_1\oplus\cdots\oplus\mathcal E_j.
\]
For $c>0$, define $\mathcal Q_{j,c}:=S_c\cap\mathcal X_j$ and
\[
\mathcal Z_{j,c}:=
\left\{u\in B_c(\widehat{\Lambda}_jc):\int_\Omega uv\,dx=0 \ \text{for every }v\in\mathcal X_{j-1}\right\}.
\]
We consider the class of odd continuous maps
\[
\Gamma_{j,c}:=\left\{h\in C\bigl(\mathcal Q_{j,c},B_c(\widehat{\Lambda}_jc)\bigr):h(-q)=-h(q)
\ \text{for every }q\in\mathcal Q_{j,c}\right\}.
\]
Since
\[
\|\nabla q\|_2^2\le\lambda_j(\Omega)c<\widehat{\Lambda}_jc\qquad\text{for every }q\in\mathcal Q_{j,c},
\]
the identity map on $\mathcal Q_{j,c}$ belongs to $\Gamma_{j,c}$. Thus, $\Gamma_{j,c}\ne\varnothing$, and we consider the following minimax problem
\[
m_{j,c}:=\inf_{h\in\Gamma_{j,c}}\max_{q\in\mathcal Q_{j,c}}E(h(q)).
\]

The following intersection property follows from the Borsuk--Ulam theorem.

\begin{lemma}\label{lem:BU}
For every $h\in\Gamma_{j,c}$, one has $h(\cQ_{j,c})\cap\cZ_{j,c}\ne\varnothing$.
\end{lemma}

\begin{proof}
Let $P_{j-1}:H_0^1(\Omega)\to\cX_{j-1}$ be the $L^2$-orthogonal projection and set $F:=P_{j-1}\circ h$. Then $F$ is continuous and odd. Since
\[
\dim\mathcal X_j-1=\dim\mathcal X_{j-1}+\dim\mathcal E_j-1\ge\dim\mathcal X_{j-1},
\]
by the Borsuk--Ulam theorem, there exists $q_h\in\mathcal Q_{j,c}$ such that $F(q_h)=0$. Hence, $\int_\Omega h(q_h)v\,dx=0$ for every $v\in\cX_{j-1}$. Since $h(q_h)\in B_c(\widehat{\Lambda}_j c)$, it follows that $h(q_h)\in\cZ_{j,c}$. Therefore, $h(\cQ_{j,c})\cap\cZ_{j,c}\ne\varnothing$.
\end{proof}

We next estimate the energy from below on $\mathcal Z_{j,c}$
and from above on $\mathcal Q_{j,c}$. Set
\[
\beta_{j,c}:=
\frac{a\lambda_j(\Omega)}{2}c
-\frac{C_{N,p}}{p}\widehat{\Lambda}_j^{q_p}c^{\frac{p}{2}},
\qquad
\alpha_{j,c}:=
\frac{a\lambda_j(\Omega)}{2}c
+\frac{b\lambda_j(\Omega)^2}{4}c^2.
\]

\begin{lemma}\label{lem:higher-sandwich}
For every fixed $j\ge2$ and every $c>0$,
\[
\inf_{u\in\cZ_{j,c}}E(u)\ge\beta_{j,c},\qquad\max_{q\in\cQ_{j,c}}E(q)\le\alpha_{j,c}.
\]
\end{lemma}

\begin{proof}
For $u\in\cZ_{j,c}$, one has $\|\nabla u\|_2^2\ge\lambda_j(\Omega)c$. Since $u\in B_c(\widehat{\Lambda}_j c)$, using \eqref{GN}, we obtain
\[
E(u)\ge\frac{a\lambda_j(\Omega)}{2}c-\frac{C_{N,p}}{p}\widehat{\Lambda}_j^{q_p}c^{\frac{p}{2}}=\beta_{j,c}.
\]
Hence, $\inf\limits_{u\in\cZ_{j,c}}E(u)\ge\beta_{j,c}$. On the other hand, for $q\in\cQ_{j,c}$, one has $\|\nabla q\|_2^2\le\lambda_j(\Omega)c$, and therefore
\[
E(q)\le\frac{a\lambda_j(\Omega)}{2}c+\frac{b\lambda_j(\Omega)^2}{4}c^2=\alpha_{j,c},
\]
which implies $\max\limits_{q\in\cQ_{j,c}}E(q)\le\alpha_{j,c}$.
\end{proof}

To verify the energy barrier condition in
Proposition~\ref{prop:localized-minimax}, we define
\[
\mathcal A_{j,c}:=\left\{u\in S_c:\Lambda_jc\le\|\nabla u\|_2^2\le
\widehat{\Lambda}_jc\right\},\qquad\delta_{j,c}:=\frac{a(\Lambda_j-\lambda_j(\Omega))}{16}c.
\]
For $j\ge2$, define
\[
c_2(j):=\min\left\{
\begin{aligned}
&1,\frac{a\bigl(\Lambda_j-\lambda_j(\Omega)\bigr)}{2b\lambda_j(\Omega)^2},\\
&\left(\frac{ap\bigl(\Lambda_j-\lambda_j(\Omega)\bigr)}{8C_{N,p}\widehat{\Lambda}_j^{q_p}}\right)^{\frac{1}{r}},\\
&\left(\frac{ap\bigl(2\lambda_j(\Omega)-\lambda_1(\Omega)-\lambda_2(\Omega)\bigr)}{4C_{N,p}\widehat{\Lambda}_j^{q_p}}\right)^{\frac{1}{r}}
\end{aligned}
\right\}.
\]

\begin{lemma}\label{lem:higher-wall}
For every fixed $j\ge2$ and every $c\in(0,c_2(j))$,
\[
\inf\limits_{u\in\mathcal A_{j,c}}E(u)>\alpha_{j,c}+4\delta_{j,c}.
\]
\end{lemma}

\begin{proof}
For $u\in\mathcal A_{j,c}$, $\Lambda_jc\le\|\nabla u\|_2^2\le\widehat{\Lambda}_jc$. Hence
\[
E(u)\ge\frac{a\Lambda_j}{2}c-\frac{C_{N,p}}{p}\widehat{\Lambda}_j^{q_p}c^{\frac{p}{2}}.
\]
Therefore,
\[
E(u)-\alpha_{j,c}-4\delta_{j,c}\ge c\left[\frac{a(\Lambda_j-\lambda_j(\Omega))}{4}-\frac{C_{N,p}}{p}\widehat{\Lambda}_j^{q_p}c^r
-\frac{b\lambda_j(\Omega)^2}{4}c\right]>0
\]
for $0<c<c_2(j)$.
\end{proof}

For $j\ge2$, define
\[
c_3(j):=\min\left\{1,\frac{a\bigl(\lambda_{j+1}(\Omega)-\lambda_j(\Omega)\bigr)}{b\lambda_j(\Omega)^2},
\left(\frac{ap\bigl(\lambda_{j+1}(\Omega)-\lambda_j(\Omega)\bigr)}{4C_{N,p}\widehat{\Lambda}_{j+1}^{q_p}}\right)^{\frac{1}{r}}\right\}.
\]

\begin{proof}[of Theorem~\ref{thm:multi}]
Fix an integer $k\ge1$. For $2\le j\le k+1$, set
\[
c_4(j):=
\begin{cases}
\displaystyle\frac{\rho^*}{\Lambda_j},
    &N=3,\ p=6,\\[1mm]
+\infty,
    &2<p<2^*.
\end{cases}
\]
With the convention that the minimum over an empty set is
$+\infty$, define
\[
c_k^*:=\min\left\{
\min_{2\le j\le k+1}
\{c_1(\Lambda_j),c_2(j),c_4(j)\},
\min_{2\le j\le k}c_3(j)
\right\}>0.
\]

Let $0<c<c_k^*$ and fix $j\in\{2,\ldots,k+1\}$.
We verify the assumptions of
Proposition~\ref{prop:localized-minimax} with $J=E$ and
\[
\mathcal G=B_c(\widehat{\Lambda}_jc),\qquad
\mathcal C=\overline{B_c(\Lambda_jc)},\qquad
Q=\mathcal Q_{j,c},\qquad
\mathcal Z=\mathcal Z_{j,c}.
\]
The identity map on $\mathcal Q_{j,c}$ belongs to
$\Gamma_{j,c}$, and Lemma~\ref{lem:BU} shows that every map
in this class meets $\mathcal Z_{j,c}$.
Together with the energy estimates of
Lemma~\ref{lem:higher-sandwich}, this verifies
{\rm(A1)--(A2)} with $\beta=\beta_{j,c}$ and
$\alpha=\alpha_{j,c}$.
Moreover, since
$\overline{\mathcal G\setminus\mathcal C}
\subset\mathcal A_{j,c}$ and $c<c_2(j)$,
the energy separation established in
Lemma~\ref{lem:higher-wall} ensures {\rm(A3)} with
$\delta=\delta_{j,c}$.
The Palais--Smale condition required in {\rm(A4)}
follows from Lemma~\ref{lem:localized-ps}.
Indeed, $\Lambda_j>\lambda_1(\Omega)$, and in the
Sobolev critical case $c<c_k^*\le c_4(j)$ gives
\[
\lambda_1(\Omega)c<\Lambda_jc<\rho^*.
\]
We therefore obtain a constrained critical point
$u_{j,c}\in\overline{B_c(\Lambda_jc)}$ such that
\[
E(u_{j,c})=m_{j,c},\qquad
\beta_{j,c}\le m_{j,c}\le\alpha_{j,c}.
\]
By the Lagrange multiplier rule, there exists
$\mu_{j,c}\in\mathbb R$ such that
$(u_{j,c},\mu_{j,c})$ solves \eqref{K}.

We next show that $u_{j,c}$ changes sign.
Since $c<c_1(\Lambda_j)$, from
Lemma~\ref{lem:one-signed-energy} it follows that every
one-signed constrained critical point in
$\overline{B_c(\Lambda_jc)}$ has energy strictly below
\[
\sigma_c:=
\frac{a}{4}\bigl(\lambda_1(\Omega)+\lambda_2(\Omega)\bigr)c.
\]
On the other hand, the choice $c<c_2(j)$ implies that
\[
\beta_{j,c}-\sigma_c
=\frac{a}{4}
\bigl(2\lambda_j(\Omega)-\lambda_1(\Omega)
      -\lambda_2(\Omega)\bigr)c
-\frac{C_{N,p}}{p}
 \widehat{\Lambda}_j^{q_p}c^{1+r}>0.
\]
Thus, $E(u_{j,c})\ge\beta_{j,c}>\sigma_c$, which excludes
the possibility that $u_{j,c}$ is one-signed. Hence, $u_{j,c}$ changes sign. 
By symmetry, $(-u_{j,c},\mu_{j,c})$ also solves \eqref{K}.

For $2\le j\le k$, from $c<c_3(j)$, we get 
\[
\alpha_{j,c}-\beta_{j+1,c}=
-\frac{a}{2}\bigl(\lambda_{j+1}(\Omega)-\lambda_j(\Omega)\bigr)c
+\frac{b\lambda_j(\Omega)^2}{4}c^2
+\frac{C_{N,p}}{p}
 \widehat{\Lambda}_{j+1}^{q_p}c^{1+r}
<0.
\]
Consequently,
\[
E(u_{j,c})\le\alpha_{j,c}
<\beta_{j+1,c}\le E(u_{j+1,c}),
\qquad 2\le j\le k.
\]
Since $E$ is even, the $k$ pairs
$\{\pm u_{2,c}\},\ldots,\{\pm u_{k+1,c}\}$
are mutually distinct.

Now we fix $j\in\{2,\ldots,k+1\}$.
Dividing the bounds for $m_{j,c}$ by $c$, we obtain
\[
\frac{a\lambda_j(\Omega)}{2}
-\frac{C_{N,p}}{p}\widehat{\Lambda}_j^{q_p}c^r
\le
\frac{E(u_{j,c})}{c}
\le
\frac{a\lambda_j(\Omega)}{2}
+\frac{b\lambda_j(\Omega)^2}{4}c.
\]
Since $r=\frac{p-2}{2}>0$, it follows that
\[
\frac{E(u_{j,c})}{c}
\longrightarrow
\frac{a\lambda_j(\Omega)}{2}
\qquad\text{as }c\to0^+.
\]
Together with
$u_{j,c}\in\overline{B_c(\Lambda_jc)}$ and
$\Lambda_j>\lambda_j(\Omega)$,
this allows us to apply
Lemma~\ref{lem:spectral-identification}
with $L=\Lambda_j$ and $\bar c=c_k^*$.
The remaining conclusions in \eqref{asym-multi} and
\eqref{profile-multi} follow.

We now prove the existence of infinitely many solutions.
If $2<p<2+\frac{8}{N}$, fix any $c>0$.
Since $q_p<2$, by the Gagliardo--Nirenberg inequality
\eqref{GN} and Young's inequality, there exists $D_c>0$
such that
\[
\frac{1}{p}\|u\|_p^p
\le
\frac{C_{N,p}}{p}
c^{\frac{p}{2}-q_p}\|\nabla u\|_2^{2q_p}
\le
\frac{b}{8}\|\nabla u\|_2^4+D_c,
\qquad u\in S_c.
\]

If $p=2+\frac{8}{N}$, set
\[
c_\infty:=
\left(\frac{b}{2C_{N,p}}\right)^{\frac{N}{4-N}}
\]
and fix $0<c<c_\infty$.
Since $q_p=2$ and $p>4$, by \eqref{GN} and the choice
of $c$, we obtain
\[
\frac{1}{p}\|u\|_p^p
\le
\frac{C_{N,p}}{p}
c^{\frac{4-N}{N}}\|\nabla u\|_2^4
\le
\frac{b}{2p}\|\nabla u\|_2^4
\le
\frac{b}{8}\|\nabla u\|_2^4,
\qquad u\in S_c.
\]
Taking $D_c=0$ in this case and using the definition
of $E$, we therefore have, in both cases,
\begin{equation}\label{eq:coercive-energy}
E(u)\ge
\frac{a}{2}\|\nabla u\|_2^2
+\frac{b}{8}\|\nabla u\|_2^4-D_c,
\qquad u\in S_c.
\end{equation}

For $j\ge2$, define
\[
\widetilde{\mathcal Z}_{j,c}
:=
\left\{
u\in S_c:
\int_\Omega uv\,dx=0
\text{ for every }v\in\mathcal X_{j-1}
\right\}
\]
and the minimax level
\[
\widetilde m_{j,c}:=
\inf_{h\in\widetilde\Gamma_{j,c}}
\max_{q\in\mathcal Q_{j,c}}E(h(q)),
\]
where 
\[\widetilde\Gamma_{j,c}:=
\left\{h\in C(\mathcal Q_{j,c},S_c):h(-q)=-h(q)\ \text{for every }q\in\mathcal Q_{j,c}\right\}.\]
By the argument in the proof of Lemma~\ref{lem:BU},
every $h\in\widetilde\Gamma_{j,c}$ satisfies
\[
h(\mathcal Q_{j,c})\cap
\widetilde{\mathcal Z}_{j,c}\ne\varnothing.
\]
For every $u\in\widetilde{\mathcal Z}_{j,c}$, we have $\|\nabla u\|_2^2\ge\lambda_j(\Omega)c$.
Combining this with \eqref{eq:coercive-energy}, we get
\[
\widetilde m_{j,c}
\ge
\inf_{u\in\widetilde{\mathcal Z}_{j,c}}E(u)
\ge
\frac{b}{8}\lambda_j(\Omega)^2c^2-D_c.
\]
Since $\lambda_j(\Omega)\to+\infty$ and $c>0$ is fixed,
we deduce that
$\widetilde m_{j,c}\to+\infty$ as $j\to\infty$.
On the other hand, the identity map belongs to
$\widetilde\Gamma_{j,c}$, and
$\|\nabla q\|_2^2\le\lambda_j(\Omega)c$
for every $q\in\mathcal Q_{j,c}$.
Using this map in the definition of
$\widetilde m_{j,c}$, we obtain
\[
\widetilde m_{j,c}
\le
\max_{q\in\mathcal Q_{j,c}}E(q)
\le
\frac{a\lambda_j(\Omega)}{2}c
+\frac{b\lambda_j(\Omega)^2}{4}c^2
=\alpha_{j,c}.
\]

For each $j\ge2$, \eqref{eq:coercive-energy} allows us to choose
$R_{j,c}>\lambda_j(\Omega)c$ such that
\[
\inf_{\substack{u\in S_c\\
                 \|\nabla u\|_2^2\ge R_{j,c}}}E(u)
>\alpha_{j,c}+4.
\]
In Proposition~\ref{prop:localized-minimax}, take
\[
\mathcal G=S_c,\qquad
\mathcal C=\overline{B_c(R_{j,c})},\qquad
Q=\mathcal Q_{j,c},\qquad
\mathcal Z=\widetilde{\mathcal Z}_{j,c},
\qquad \delta=1.
\]
Using the intersection property and the lower and upper
energy bounds, we verify {\rm(A1)--(A2)} with
\[
\beta=\frac{b}{8}\lambda_j(\Omega)^2c^2-D_c,
\qquad
\alpha=\alpha_{j,c}.
\]
By the definition of $\mathcal C$ and our choice of
$R_{j,c}$, we have {\rm(A3)}.
Since $p<2^*$ in both cases, we conclude from
Lemma~\ref{lem:localized-ps} that {\rm(A4)} holds.
We therefore apply
Proposition~\ref{prop:localized-minimax} to obtain
a constrained critical point
$\widetilde u_{j,c}\in\overline{B_c(R_{j,c})}$ satisfying
$E(\widetilde u_{j,c})=\widetilde m_{j,c}$.

To show that $\widetilde u_{j,c}$ changes sign for all
sufficiently large $j$, we estimate the energies of
one-signed normalized solutions of \eqref{K} with mass $c$.
Let $u\in S_c$ be any one-signed solution of \eqref{K}
with Lagrange multiplier $\mu$, and set
$A:=\|\nabla u\|_2^2$.
Replacing $u$ by $-u$ if necessary, we may assume that
$u\ge0$ in $\Omega$.
Since $u\not\equiv0$, the strong maximum principle implies that $u>0$ in $\Omega$.
Testing \eqref{K} with $\phi_1$, as in the proof of
Lemma~\ref{lem:one-signed-energy}, we obtain
\[
\mu<(a+bA)\lambda_1(\Omega).
\]
Testing \eqref{K} with $u$, we also have
\[
(a+bA)A=\|u\|_p^p+\mu c.
\]
Combining these two relations and using \eqref{GN},
we find that
\[
(a+bA)\bigl(A-\lambda_1(\Omega)c\bigr)
<
\|u\|_p^p
\le
C_{N,p}c^{\frac{p}{2}-q_p}A^{q_p}.
\]

We first consider the case $q_p<2$.
If $A\ge2\lambda_1(\Omega)c$, then
\[
\frac{b}{2}A^2
\le
(a+bA)\bigl(A-\lambda_1(\Omega)c\bigr)
<
C_{N,p}c^{\frac{p}{2}-q_p}A^{q_p}.
\]
Since $2-q_p>0$, it follows that
\[
A<
\left(
\frac{2C_{N,p}}{b} c^{\frac{p}{2}-q_p}
\right)^{\frac{1}{2-q_p}}.
\]
Combining this estimate with the alternative
$A<2\lambda_1(\Omega)c$, we conclude
that every one-signed constrained critical point satisfies
\[
A\le
\max\left\{
2\lambda_1(\Omega)c,\,
\left(
\frac{2C_{N,p}}{b} c^{\frac{p}{2}-q_p}
\right)^{\frac{1}{2-q_p}}
\right\}.
\]

We next consider the case $q_p=2$.
By the Poincar\'e inequality,
$A\ge\lambda_1(\Omega)c$.
Using this bound and $c<c_\infty$, we obtain
\[
bA\bigl(A-\lambda_1(\Omega)c\bigr)
\le
(a+bA)\bigl(A-\lambda_1(\Omega)c\bigr)
<
C_{N,p}c^{\frac{4-N}{N}}A^2
\le\frac{b}{2}A^2.
\]
Dividing by $bA>0$, we find
$A-\lambda_1(\Omega)c<\frac{A}{2}$,
and hence $A<2\lambda_1(\Omega)c$.

In both cases, the bounds on $A$ and
\[
E(u)\le\frac{a}{2}A+\frac{b}{4}A^2
\]
show that the energies of all one-signed constrained
critical points in $S_c$ are bounded above.
Since $E(\widetilde u_{j,c})=\widetilde m_{j,c}\to+\infty$,
the points $\widetilde u_{j,c}$ change sign for all
sufficiently large $j$.
Selecting a subsequence with strictly increasing energies
and using the evenness of $E$, we obtain infinitely many
mutually distinct pairs.
This completes the proof.
\end{proof}

\section{Least energy sign-changing normalized solutions}
\label{glne:section}

This section is devoted to the proof of
Theorem~\ref{glne:main}.
Recall that $d_c$ is the infimum of $E$ over
$\mathcal K_{c}^{\mathrm{nod}}$, which is nonempty for
sufficiently small $c>0$ by Theorem~\ref{thm:multi}.
We begin with an a priori energy estimate and a
compactness result that preserves the sign-changing property.

For $2<p\le2+\frac{8}{N}$, define
\[
\tau_p:=\min\left\{\frac{p-2}{2},\,\frac{p}{2}-q_p\right\}>0
\]
and
\[
c_5:=
\min\left\{1,
\left[
\frac{p}{C_{N,p}\bigl(1+\lambda_1(\Omega)^{q_p-1}\bigr)}
\min\left\{\frac{a}{4},\frac{b}{8}\right\}
\right]^{\frac{1}{\tau_p}}
\right\}.
\]
In the remaining cases of Theorem~\ref{glne:main},
set $c_5=1$.

\begin{lemma}\label{glne:energy-bound}
Under any of the assumptions of Theorem~\ref{glne:main}, for every $0<c<c_5$ and every constrained critical point $u\in S_c$, $E(u)\ge\frac{a}{4} \|\nabla u\|^2_2$. Moreover, for every $M>0$, if $E(u)\le Mc$, then $\|\nabla u\|^2_2\le\frac{4M}{a}c$.
\end{lemma}

\begin{proof}
We first consider the case $2<p<2+\frac{8}{N}$. Then, $q_p\in(0,2)$.

If $q_p\in(0,1]$, by the Gagliardo--Nirenberg inequality \eqref{GN} and the Poincar\'e inequality, we have
\[
\|u\|_p^p\le C_{N,p}c^{\frac{p}{2}-q_p}\|\nabla u\|_2^{2q_p}\le C_{N,p}\lambda_1(\Omega)^{q_p-1}c^{\frac{p-2}{2}}\|\nabla u\|_2^2.
\]
It follows that
\[
E(u)\ge\left[\frac{a}{2}-\frac{C_{N,p}}{p}\lambda_1(\Omega)^{q_p-1}c^{\frac{p-2}{2}}\right]\|\nabla u\|_2^2
+\frac{b}{4}\|\nabla u\|_2^4.
\]
By the choice of $c_5$, we obtain
\[
E(u)\ge\frac{a}{4}\|\nabla u\|_2^2+\frac{b}{4}\|\nabla u\|_2^4\ge\frac{a}{4}\|\nabla u\|_2^2.
\]

If $q_p\in(1,2)$, using $t^{q_p}\le t+t^2$ for $t\ge0$, we have
\[
\|u\|_p^p\le C_{N,p}c^{\frac{p}{2}-q_p}\left(\|\nabla u\|_2^2+\|\nabla u\|_2^4\right).
\]
Consequently,
\[
E(u)\ge\left(\frac{a}{2}-\frac{C_{N,p}}{p}c^{\frac{p}{2}-q_p}\right)\|\nabla u\|_2^2
+\left(\frac{b}{4}-\frac{C_{N,p}}{p}c^{\frac{p}{2}-q_p}\right)\|\nabla u\|_2^4.
\]
By the choice of $c_5$, we obtain
\[
E(u)\ge\frac{a}{4}\|\nabla u\|_2^2+\frac{b}{8}\|\nabla u\|_2^4\ge\frac{a}{4}\|\nabla u\|_2^2.
\]

We next consider the case $p=2+\frac{8}{N}$. Then $q_p=2$ and $\tau_p=\frac{4}{N}-1>0$. By \eqref{GN} and the choice of $c_5$, we have
\[
E(u)\ge\frac{a}{2}\|\nabla u\|_2^2+
\left(\frac{b}{4}-\frac{C_{N,p}}{p}c^{\frac{4}{N}-1}\right)\|\nabla u\|_2^4
\ge\frac{a}{2}\|\nabla u\|_2^2+\frac{b}{8}\|\nabla u\|_2^4\ge\frac{a}{4}\|\nabla u\|_2^2.
\]

Finally, suppose that $\Omega$ is star-shaped and that
either $2+\frac{8}{N}<p<2^*$ or $(N,p)=(3,6)$.
After a translation, we may assume that $x\cdot\nu\ge0$
on $\partial\Omega$, where $\nu$ denotes the outward
unit normal.
Let $\mu$ be the Lagrange multiplier associated with $u$,
and recall that $M_u=a+b\|\nabla u\|_2^2$.
Since $M_u$ is independent of $x$, the regularity established
above allows us to multiply the equation by $x\cdot\nabla u$
and integrate by parts. This calculation is valid for
$1\le N\le3$ and gives the Pohozaev identity
(see \cite[Theorem~B.1]{Willem-1996}):
\[
\frac{N-2}{2}M_u\|\nabla u\|_2^2
-\frac{N}{p}\|u\|_p^p-\frac{N\mu}{2}c
+\frac{M_u}{2}
\int_{\partial\Omega}
(x\cdot\nu)|\partial_\nu u|^2\,d\sigma
=0.
\]
Testing \eqref{K} with $u$, we also have
\[
M_u\|\nabla u\|_2^2=\|u\|_p^p+\mu c.
\]
Combining these identities and using $x\cdot\nu\ge0$ on $\partial\Omega$, we obtain
\[
M_u\|\nabla u\|_2^2-\frac{N(p-2)}{2p}\|u\|_p^p
=
\frac{M_u}{2}
\int_{\partial\Omega}(x\cdot\nu)|\partial_\nu u|^2\,d\sigma
\ge0.
\]
Consequently,
\[
\frac{1}{p}\|u\|_p^p
\le
\frac{2}{N(p-2)}
\left(a\|\nabla u\|_2^2+b\|\nabla u\|_2^4\right),
\]
and hence
\[
E(u)\ge
a\left(\frac{1}{2}-\frac{2}{N(p-2)}\right)\|\nabla u\|_2^2
+
b\left(\frac{1}{4}-\frac{2}{N(p-2)}\right)\|\nabla u\|_2^4.
\]
Since $N(p-2)>8$ in both cases, it follows that
$E(u)\ge\frac{a}{4}\|\nabla u\|_2^2$.

We have therefore established this lower bound in all the cases under consideration. If in addition \(E(u)\le Mc\), then
\[
\|\nabla u\|_2^2\le\frac{4}{a}E(u)\le\frac{4M}{a}c.
\]
\end{proof}

We next establish a compactness result that preserves
the sign-changing property.
\begin{lemma}\label{glne:compactness}
Let $c>0$ and $\rho>\lambda_1(\Omega)c$, and let
$\{u_n\}\subset
\mathcal K_{c}^{\mathrm{nod}}\cap\overline{B_c(\rho)}$.
Assume either $2<p<2^*$, or $N=3$, $p=6$, and
$\rho\le\rho^*$.
Then, after passing to a subsequence,
$u_n\to u$ in $H_0^1(\Omega)$ for some
$u\in
\mathcal K_{c}^{\mathrm{nod}}\cap\overline{B_c(\rho)}$.
\end{lemma}

\begin{proof}
Since $\{u_n\}$ is bounded in $H_0^1(\Omega)$,
the sequence $\{E(u_n)\}$ is bounded.
After passing to a subsequence, we may assume that
$E(u_n)\to d$ for some $d\in\mathbb R$.
As each $u_n$ is a constrained critical point, by
Lemma~\ref{lem:localized-ps} we get
\[
u_n\to u
\quad\text{in }H_0^1(\Omega)
\]
for some $u\in\overline{B_c(\rho)}$.
By the Sobolev inequality \eqref{Sobolev}, we also have
\[
\|u_n-u\|_p
\le\mathcal S_p\|u_n-u\|\to0,
\]
which covers the case $N=3$ and $p=6$.

Let $\mu_n$ be the multiplier associated with $u_n$.
Then
\begin{equation}\label{glne:un-equation}
-(a+b\|\nabla u_n\|_2^2)\Delta u_n
=
|u_n|^{p-2}u_n+\mu_nu_n
\quad\text{in }\Omega.
\end{equation}
Testing \eqref{glne:un-equation} with $u_n$, we obtain
\[
\mu_n
=
\frac{(a+b\|\nabla u_n\|_2^2)\|\nabla u_n\|_2^2-\|u_n\|_p^p}{c}.
\]
Using the strong convergence established above it follows that
$\mu_n\to\mu$ for some $\mu\in\mathbb R$.
In particular, there exists $C>0$ such that
$|\mu_n|\le C$ for every $n$.
Passing to the limit in the weak formulation of
\eqref{glne:un-equation}, we find that $u$ is a
constrained critical point of $E|_{S_c}$.

It remains to show that $u$ changes sign.
Suppose, by contradiction, that $u\ge0$.
Since $\|u\|_2^2=c>0$, by the strong maximum principle we obtain
 $u>0$ in $\Omega$.
As each $u_n$ changes sign, we can define
\[
z_n:=\frac{u_n^-}{\|u_n^-\|_2}.
\]
Then, $z_n\ge0$ and $\|z_n\|_2=1$.
Testing \eqref{glne:un-equation} with $-u_n^-$
and dividing by $\|u_n^-\|_2^2$, we obtain
\begin{equation}\label{eq:normalized-negative-part}
(a+b\|\nabla u_n\|_2^2)\|\nabla z_n\|_2^2
=
\int_\Omega |u_n|^{p-2}z_n^2\,dx+\mu_n.
\end{equation}

We next show that $\{z_n\}$ is bounded in $H_0^1(\Omega)$.
Since $u_n\to u$ in $L^p(\Omega)$ and $u\ge0$, we have
\[
0\le u_n^-\le |u_n-u|
\quad\text{a.e. in }\Omega,
\qquad
\|u_n^-\|_p\le\|u_n-u\|_p\to0.
\]
By the definition of $z_n$, H\"older's inequality,
and \eqref{Sobolev}, it follows that 
\begin{equation}\label{glne:negative-part-estimate}
\begin{aligned}
\int_\Omega |u_n|^{p-2}z_n^2\,dx
&=
\int_\Omega (u_n^-)^{p-2}z_n^2\,dx\\
&\le
\|u_n^-\|_p^{p-2}\|z_n\|_p^2\\
&\le
\mathcal S_p^2\|u_n^-\|_p^{p-2}
\|\nabla z_n\|_2^2\\
&\le
\frac{a}{2}\|\nabla z_n\|_2^2
\end{aligned}
\end{equation}
for all sufficiently large $n$, since $p>2$ and
$\|u_n^-\|_p\to0$.
Combining \eqref{eq:normalized-negative-part} with
\eqref{glne:negative-part-estimate}, and using
$a+b\|\nabla u_n\|_2^2\ge a$ and $\mu_n\to\mu$,
we obtain
\begin{equation}\label{glne:zn-gradient}
\frac{a}{2}\|\nabla z_n\|_2^2
\le \mu_n\le |\mu|+1
\end{equation}
for all sufficiently large $n$.
By \eqref{glne:zn-gradient} and the compact embedding
$H_0^1(\Omega)\hookrightarrow L^2(\Omega)$, we may assume,
after passing to a further subsequence, that
\[
z_n\to z\quad\text{in }L^2(\Omega),
\qquad z\ge0,\quad \|z\|_2=1.
\]

Since $u_n\to u$ in $L^2(\Omega)$ and $u\ge0$,
we have $\|u_n^-\|_2\to0$.
By the definition of $z_n$,
\[
\int_\Omega u_nz_n\,dx
=-\|u_n^-\|_2\to0.
\]
Using the strong $L^2$ convergence of both $u_n$ and $z_n$,
we may pass to the limit and obtain
$\int_\Omega uz\,dx=0$.
This contradicts $u>0$ in $\Omega$, $z\ge0$,
and $\|z\|_2=1$.
Thus $u^-\not\equiv0$.
Applying the same argument to $\{-u_n\}$ gives
$u^+\not\equiv0$.
Therefore,
$u\in\mathcal K_c^{\mathrm{nod}}
\cap\overline{B_c(\rho)}$.
This completes the proof.
\end{proof}

\begin{proof}[of Theorem~\ref{glne:main}]
We first show that $d_c$ is attained in case~{\rm(i)}. 
Fix $c>0$ if $2<p<2+\frac{8}{N}$, and fix
$0<c<c_\infty$ if $p=2+\frac{8}{N}$.
By Theorem~\ref{thm:multi},
$\mathcal K_{c}^{\mathrm{nod}}\ne\varnothing$.
For these values of $c$, \eqref{eq:coercive-energy} holds with some $D_c\ge0$:
\[
E(u)\ge
\frac{a}{2}\|\nabla u\|_2^2
+\frac{b}{8}\|\nabla u\|_2^4-D_c,
\qquad u\in S_c.
\]
Thus $-\infty<d_c<+\infty$. Let $\{u_n\}\subset\mathcal K_{c}^{\mathrm{nod}}$
be a minimizing sequence.
By \eqref{eq:coercive-energy}, $\{u_n\}$ is bounded
in $H_0^1(\Omega)$. Choose $\rho>\lambda_1(\Omega)c$ large enough that
$\{u_n\}\subset\overline{B_c(\rho)}$.
Since $p<2^*$, Lemma~\ref{glne:compactness}
applies with this $\rho$.
After passing to a subsequence, we have
$u_n\to w_c$ in $H_0^1(\Omega)$ for some $w_c\in\mathcal K_{c}^{\mathrm{nod}}$. By the continuity of $E$,
\[
E(w_c)=\lim_{n\to\infty}E(u_n)=d_c.
\]

Let $c_1^*>0$ be the threshold given by
Theorem~\ref{thm:multi} with $k=1$.
We next consider cases~{\rm(ii)} and~{\rm(iii)}.
For every $0<c<c_1^*$, the solution constructed in its proof
satisfies $u_{2,c}\in\mathcal K_{c}^{\mathrm{nod}}
\cap\overline{B_c(\Lambda_2c)}$.
Since $\|\nabla u_{2,c}\|_2^2\le\Lambda_2c$, by the definition
of $E$, we get
\[
d_c\le E(u_{2,c})
\le\frac{a\Lambda_2}{2}c
+\frac{b\Lambda_2^2}{4}c^2
\le C_0c,
\qquad 0<c<c_1^*,
\]
where
\[
C_0:=\frac{a\Lambda_2}{2}
+\frac{b\Lambda_2^2c_1^*}{4}
\]
is independent of $c$. Set $M:=C_0+1$ and $L:=\frac{4M}{a}$, so that
$L>2\Lambda_2>\lambda_1(\Omega)$.
Choose $0<c^*<\min\{c_1^*,c_5\}$ and require
$Lc^*<\rho^*$ when $N=3$ and $p=6$.
Fix $0<c<c^*$.
Then $\rho=Lc$ satisfies $\rho>\lambda_1(\Omega)c$,
and also $\rho<\rho^*$ in the Sobolev critical case.
For every $u\in\mathcal K^{\mathrm{nod}}_{c}$, using Lemma~\ref{glne:energy-bound} and the Poincar\'e inequality, we have
\[
E(u)
\ge
\frac{a}{4}\|\nabla u\|_2^2
\ge
\frac{a\lambda_1(\Omega)}{4}c.
\]
It follows that
\[
\frac{a\lambda_1(\Omega)}{4}c
\le d_c\le C_0c.
\]
Choose a minimizing sequence $\{u_n\}\subset\mathcal K^{\mathrm{nod}}_{c}$ such that
\[
d_c\le E(u_n)\le d_c+\frac{c}{n}\le Mc.
\]
Since each $u_n$ is a constrained critical point, in view of Lemma~\ref{glne:energy-bound}, we get
$\|\nabla u_n\|_2^2\le Lc$.
By Lemma~\ref{glne:compactness} with $\rho=Lc$, after
passing to a subsequence, we have
$u_n\to w_c$ in $H_0^1(\Omega)$
for some $w_c\in\mathcal K_{c}^{\mathrm{nod}}$.
By the continuity of $E$,
\[
E(w_c)=\lim_{n\to\infty}E(u_n)=d_c.
\]

It remains to prove \eqref{asym-least-nodal}
and \eqref{profile-least-nodal} in all three cases.
For the asymptotic analysis in case~{\rm(i)}, we also
define $C_0$, $M$, and $L$ by the same formulas as above
and choose an auxiliary $c^*>0$ such that
$c^*<\min\{c_1^*,c_5\}$.
When $p=2+\frac{8}{N}$, we additionally require
$c^*<c_\infty$.
For $0<c<c^*$, the solution $u_{2,c}$ exists in
case~{\rm(i)} and satisfies
$d_c\le E(u_{2,c})\le C_0c$ by the same estimate as above.
We may therefore restrict to $0<c<c^*$ in all three cases
when proving the limits as $c\to0^+$.
Let $w_c$ be any least energy sign-changing normalized
solution, and let $\mu_c$ be its associated Lagrange multiplier. By Lemma~\ref{glne:energy-bound}, we have
\[
\|\nabla w_c\|_2^2
\le\frac{4}{a}d_c
\le\frac{4C_0}{a}c<Lc.
\]
Set $v_c:=\frac{w_c}{\sqrt c}$.
Then $\|v_c\|_2=1$ and $\|\nabla v_c\|_2^2<L$.
By \eqref{Sobolev}, $\|v_c\|_p$ is uniformly bounded
as $c\to0^+$.
Using the definition of $E$ and testing \eqref{K}
with $w_c$, we obtain
\begin{equation}\label{eq:least-nodal-expansion}
\frac{d_c}{c}
=\frac{a}{2}\|\nabla v_c\|_2^2+o(1),
\qquad
\mu_c=a\|\nabla v_c\|_2^2+o(1).
\end{equation}
In view of $d_c\le E(u_{2,c})$, \eqref{asym-multi} implies that
\begin{equation}\label{eq:least-nodal-limsup}
\limsup_{c\to0^+}\frac{d_c}{c}
\le\frac{a\lambda_2(\Omega)}{2}.
\end{equation}
Let $c_n\to0^+$ be arbitrary. Since $\{v_{c_n}\}$ is
bounded in $H_0^1(\Omega)$, after passing to a subsequence
we have
\[
\|\nabla v_{c_n}\|_2^2\to\theta,
\qquad
v_{c_n}\rightharpoonup v\quad\text{in }H_0^1(\Omega),\qquad
v_{c_n}\to v\quad\text{in }L^2(\Omega).
\]
Using $\|v_{c_n}\|_2=1$, the Poincar\'e inequality,
\eqref{eq:least-nodal-expansion}, and
\eqref{eq:least-nodal-limsup}, we obtain
\[
\|v\|_2=1,
\qquad
\lambda_1(\Omega)\le\theta\le\lambda_2(\Omega),
\qquad
\mu_{c_n}\to a\theta.
\]
For each $n$, $v_{c_n}$ satisfies
\begin{equation}\label{eq:vn}
-\left(a+bc_n\|\nabla v_{c_n}\|_2^2\right)\Delta v_{c_n}
=c_n^r|v_{c_n}|^{p-2}v_{c_n}+\mu_{c_n}v_{c_n}.
\end{equation}
By H\"older's inequality and \eqref{Sobolev}, for every
$\psi\in H_0^1(\Omega)$,
\[
\left|c_n^r\int_\Omega
|v_{c_n}|^{p-2}v_{c_n}\psi\,dx\right|
\le c_n^r\mathcal S_p\|v_{c_n}\|_p^{p-1}\|\psi\|\to0.
\]
Moreover,
$bc_n\|\nabla v_{c_n}\|_2^2\to0$.
Testing \eqref{eq:vn} with
$\psi$ and letting $n\to\infty$, we obtain
\[
a\int_\Omega\nabla v\cdot\nabla\psi\,dx
=a\theta\int_\Omega v\psi\,dx
\qquad\text{for every }\psi\in H_0^1(\Omega).
\]
Hence, $-\Delta v=\theta v$ in $\Omega$. 
Testing this equation with $v$, we obtain
$\|\nabla v\|_2^2=\theta$. Combining $v_{c_n}\rightharpoonup v$ in $H_0^1(\Omega)$
with $\|\nabla v_{c_n}\|_2^2\to\|\nabla v\|_2^2$,
we obtain $v_{c_n}\to v$ in $H_0^1(\Omega)$.
Since $\theta$ is an eigenvalue in
$[\lambda_1(\Omega),\lambda_2(\Omega)]$ and
$\lambda_1(\Omega),\lambda_2(\Omega)$ are the first two
distinct Dirichlet eigenvalues, we conclude that
$\theta\in\{\lambda_1(\Omega),\lambda_2(\Omega)\}$.

We now show that $\theta=\lambda_2(\Omega)$.
Suppose, for contradiction, that
$\theta=\lambda_1(\Omega)$.
Since the first eigenvalue is simple, after changing signs
if necessary, we may assume that $v=\phi_1>0$ in $\Omega$.
As each $w_{c_n}$ changes sign, we can define
\[
z_n:=\frac{w_{c_n}^-}{\|w_{c_n}^-\|_2},
\qquad
z_n\ge0,\quad \|z_n\|_2=1.
\]
Testing the equation for $w_{c_n}$ with $-w_{c_n}^-$
and dividing by $\|w_{c_n}^-\|_2^2$, we obtain
\begin{equation}\label{eq:least-nodal-negative-part}
\left(a+b\|\nabla w_{c_n}\|_2^2\right)
\|\nabla z_n\|_2^2
=
\int_\Omega |w_{c_n}|^{p-2}z_n^2\,dx
+\mu_{c_n}.
\end{equation}
By H\"older's inequality, the Sobolev inequality, and
$\|\nabla w_{c_n}\|_2^2\le Lc_n$, we have
\begin{equation}\label{eq:zn}
\begin{aligned}
\int_\Omega |w_{c_n}|^{p-2}z_n^2\,dx
&\le \|w_{c_n}\|_p^{p-2}\|z_n\|_p^2\\
&\le \mathcal S_p^p
      \|\nabla w_{c_n}\|_2^{p-2}\|\nabla z_n\|_2^2\\
&\le \mathcal S_p^p(Lc_n)^{\frac{p-2}{2}}
      \|\nabla z_n\|_2^2.
\end{aligned}
\end{equation}
Since $p>2$ and $c_n\to0^+$,
\[
\mathcal S_p^p(Lc_n)^{\frac{p-2}{2}}\le\frac{a}{2}
\]
for all sufficiently large $n$. Combining \eqref{eq:least-nodal-negative-part}--\eqref{eq:zn}
with $\mu_{c_n}\to a\lambda_1(\Omega)$, we obtain, for all sufficiently large $n$,
\[
\|\nabla z_n\|_2^2
\le\frac{2}{a}\mu_{c_n}
\le2\lambda_1(\Omega)+\frac{2}{a}.
\]
By the compact embedding
$H_0^1(\Omega)\hookrightarrow L^2(\Omega)$,
we may therefore assume, after passing to a subsequence, that
\[
z_n\to z\quad\text{in }L^2(\Omega),
\qquad
z\ge0,\quad \|z\|_2=1.
\]

On the other hand, since $v_{c_n}\to\phi_1$ in $L^2(\Omega)$
and $\phi_1>0$ in $\Omega$, we have
$\|v_{c_n}^-\|_2\to0$.
By the definitions of $v_{c_n}$ and $z_n$, we get
\[
\int_\Omega v_{c_n}z_n\,dx
=-\|v_{c_n}^-\|_2\to0.
\]
Using $v_{c_n}\to\phi_1$ and $z_n\to z$ in
$L^2(\Omega)$, we obtain
$\int_\Omega\phi_1z\,dx=0$.
This contradicts $\phi_1>0$ in $\Omega$, $z\ge0$, and
$\|z\|_2=1$.
Hence, $\theta=\lambda_2(\Omega)$ and
$v\in\mathbb S(\cE_2)$.
Since $c_n\to0^+$ was arbitrary,
\eqref{asym-least-nodal} and \eqref{profile-least-nodal}
follow.
\end{proof}

\section*{Acknowledgements}
This work is partially supported by the NSFC (12471102), NSF of Jilin Province (20250102004JC), and Research Project of the Education Department of Jilin Province (JJKH20250296KJ).

\section*{Data availability}
Data sharing is not applicable to this article as no datasets were generated or analysed during the current study.

\section*{Declarations}
Conflict of interest. The authors declare that they have no conflict of interest.

\end{document}